\input amstex

\loadeufm
\loadmsbm
\loadeufm

\documentstyle{amsppt}
\input amstex
\catcode `\@=11
\def\logo@{}
\catcode `\@=12
\magnification \magstep1
\NoRunningHeads
\NoBlackBoxes
\TagsOnLeft

\def \={\ = \ }
\def \+{\ +\ }
\def \-{\ - \ }

\def \b|{\big |}

\def \g1{\Gamma_1}

\def \nfp{\demo\nofrills{Proof:\usualspace\usualspace }}

\def \stab{\text{stab }}

\def\rarr#1#2{\smash{\mathop{\hbox to .5in{\rightarrowfill}}
 	 \limits^{\scriptstyle#1}_{\scriptstyle#2}}}

\def\larr#1#2{\smash{\mathop{\hbox to .5in{\leftarrowfill}}
	  \limits^{\scriptstyle#1}_{\scriptstyle#2}}}

\def\swarr#1#2 {\llap{$\scriptstyle #1$}  \swarrow
  	\vcenter to .5in{}\rlap{$\scriptstyle #2$}}

\topmatter
\title Makanin-Razborov Diagrams over Free Products 
\endtitle
\author
\centerline{ 
E. Jaligot${}^{1,2}$ and Z. Sela${}^{3,4}$}
\endauthor
\footnote""{${}^1$Universite de Lyon - CNRS and Universite Lyon 1.}
\footnote""{${}^2$Partially supported by a French ANR JC05-47038.} 
\footnote""{${}^3$Hebrew University, Jerusalem 91904, Israel.}
\footnote""{${}^4$Partially supported by an Israel academy of sciences fellowship.} 
\abstract\nofrills{}
This paper is the first in a sequence on the first order theory of free
products. 
In the first paper we generalize the analysis of systems of equations
over free and (torsion-free) hyperbolic groups, and analyze systems of equations
over free products. To do that we introduce limit groups over the class
of free products, and show that a finitely presented group has a canonical
(finite) collection of maximal limit quotients. We further extend this finite
collection and associate a Makanin-Razborov diagram over free products with
every f.p.\ group. This MR diagram encodes all the quotients of a given
f.p.\ group that are free products, all its homomorphisms into free
products, and equivalently all the solutions to a given system of equations
over a free product.
\endabstract
\endtopmatter

\document

\baselineskip 12pt

Sets of solutions to equations defined over a free group have been studied
extensively.  
Considerable progress in the study of
such sets of solutions was made by G. S. Makanin, who constructed
an algorithm that decides if a system of equations defined over a free group
has a solution [Ma], and showed that the universal and positive theories
of a free group are decidable. A. A. Razborov
was able to give a description of the entire set of solutions to a system of
equations defined over a free group [Ra2], a description that was further developed
by O. Kharlampovich and A. Myasnikov [Kh-My]. 

In [Se1] a geometric approach to the study of sets of solutions to systems of equations over a
free group is presented. This was generalized in [Se3] for systems of equations over (torsion-free)
hyperbolic groups, in [Al] to systems of equations over limit groups, and in [Gr] to systems of equations
over toral relatively hyperbolic groups.

In this paper we generalize part of the techniques and results that were obtained over free
groups to study systems of equations over arbitrary free products. Let $\Sigma$ be a system of
equations which is defined over a free product, $A*B$:
$$\align
w_1(x_1,\ldots,x_n) = 1 \\
& \vdots \hfill \\
w_s(x_1,\ldots,x_n) = 1 \\
\endalign$$
Following [Ra1] we set the associated f.p.\ group $G(\Sigma)$ to be:
$$
G(\Sigma) \ = \ < \, x_1,\ldots,x_n \, | \, w_1,\ldots,w_s \, >
$$

\noindent
Clearly, every solution of the system $\Sigma$ corresponds to a homomorphism
$h:G(\Sigma) \to A*B$,  and every such homomorphism 
corresponds to a solution of the system $\Sigma$. Therefore, the study of 
sets of solutions to systems of equations over the free product $A*B$ is equivalent to
the study of all the homomorphisms from a fixed f.p.\ group $G$ into $A*B$.

We further generalize our point of view, and instead of the set of homomorphisms from a 
given f.p.\ group $G(\Sigma)$ into a particular free product, we study the set of all
the homomorphisms from the f.p.\ group $G(\Sigma)$ into all possible free products. By
Kurosh subgroup theorem, this 
is equivalent to the study of all the quotients of a given f.p.\ group, $G(\Sigma)$, that
are free products. 

To analyze the set of free product quotients of a given f.p.\ group, we generalize the notion of
limit groups (over free groups), and define $limit$ $groups$ $over$ $free$ $products$. The definition
over free products (definition 1) is a generalization of the definition of limit groups over free groups,
but  with each limit group over free products, $L$, there is an additional structure, a subset of
conjugacy classes in the limit group $L$, that are called $elliptics$, that are forced to be mapped to
conjugates of the factors in any homomorphism from the limit group into a free product.

After proving some basic properties of limit groups over free products, we associate with them a canonical
virtually abelian JSJ decomposition (theorem 11). Limit groups over free products do not
satisfy the d.c.c.\ that hold for limit groups over free and hyperbolic groups. Still, in theorem 13 we prove a
basic d.c.c.\ that holds for such limit groups, and applies to descending chain of limit groups over free
products, in which the maps between successive limit groups are proper epimorphisms that do not map  non-trivial
elliptic elements to the identity element. 

This d.c.c.\  allows us to associate a $resolution$  with each limit group over free
products (theorem 18). We further
define a natural partial order on the set of limit quotients over free products of a given f.p.\ group, and
prove that there are finitely many (equivalence classes of) maximal limit quotients (over
free products) of a f.p.\ group. Finally we extend each of the maximal limit quotients with finitely many
resolutions and obtain a Makanin-Razborov diagram of a f.p.\ group over free products.

The diagram that we associate with a f.p.\ group encodes all the quotients of the given f.p.\ group that
are free products. Unfortunately, our construction is not canonical, and we state a natural conjecture that
if answered affirmatively will enable one to construct a canonical diagram. Also, the construction uses the
finite presentability of the group in question in an essential way. Hence, encoding the set of free product 
quotients of a f.g.\ group is left  open.

The Makanin-Razborov diagram over free products is the first step towards the analysis of the
first order theory of free products that will appear in the  sequel. This  study  was motivated by a question 
of the first author on the stability of a free product of stable groups. We expect that some of the notions
and constructions that appear in this paper (and in the sequel) can be generalized to other classes of
groups, e.g. acylindrical splittings of f.p.\ groups, and various classes of relatively hyperbolic groups.

\vglue 1.5pc
\centerline{\bf{\S1. Limit Groups over Free Products}}
\medskip

We start the analysis of systems of equations over free products 
with the definition of a limit group over the set of free products. The definition generalizes
the corresponding ones for free, hyperbolic, and relatively hyperbolic groups, but it associates
with a limit group an additional structure - it's collection of conjugacy classes of elliptic elements.
Also, note that unlike the case of a free or a hyperbolic group, we consider limit groups over the entire
class of free products, and not necessarily over a given one.

\vglue 1pc
\proclaim{Definition 1} 
Let $\{A_n\}$ and $\{B_n\}$ be two sequences of groups (not necessarily finitely
generated), and let $G$ be a finitely generated group. We say that a sequence of
homomorphisms, $\{h_n:G \to A_n*B_n\}$, is a convergent sequence, if the following 
conditions hold:
\roster
\item"{(i)}" for each $g \in G$ there exists some index $n_g>0$, so that for every $n>n_g$, 
$h_n(g)=1$, or for every $n>n_g$, $h_n(g) \neq 1$.  

\item"{(ii)}" for each $g \in G$ there exists some index $n^e_g>0$, so that for every $n>n^e_g$, 
$h_n(g)$ is elliptic in the free product $A_n*B_n$ (i.e., it is contained in a conjugate of $A_n$ or
$B_n$), or for every $n>n^e_g$, $h_n(g)$  is not elliptic in $A_n*B_n$. 
\endroster

With the convergent sequence we associate its $stable$ $kernel$ that is defined to be:
$$
K=\{g \in G \, | \, \exists n_g \, \forall n>n_g \, h_n(g)=1  \}
$$
and the associate limit group:
 $L=G/K$, which we call a limit group over (the collection of) free products, 
and set $\eta : G \to L$ to be the
 natural quotient map.

With the limit group $L$ we associate an additional structure, its  collection of conjugacy classes that
are stably elliptic, i.e.:
$$
E_L=\{\ell \in L \, | \, \exists g \in G \, \eta(g) = \ell  \, \exists n_g>0 \, \forall n>n_g \, h_n(g) \ is \
elliptic \}
$$
Note that by definition if $\eta(g_1)=\eta(g_2)$, then $g_1$ is stably elliptic iff $g_2$ is stably elliptic.
Also, note that every f.g.\ group can be a limit group over free products, as given a finitely generated group 
$G$, we can look at the free product $G*B$, for some non-trivial group $B$,
with the fixed sequence of homomorphisms that map $G$ identically
onto $G$ in the free product $G*B$. Note that in this tautological case, the entire (limit) group $G$ 
is set to be elliptic.
\endproclaim

Given a convergent sequence of homomorphisms one can  pass to a 
subsequence that converges into
a (possibly trivial) action of the associated limit group on some real tree. 

\noindent
Let $A$ and $B$ be non-trivial groups (not necessarily finitely generated). With the free product, $A*B$,
we can naturally associate its Bass-Serre tree. Let $G$ be  a f.g.\ group $G=<g_1,\ldots,g_m>$,  let
$\{A_n,B_n\}$ be a sequence
of pairs of non-trivial groups, and let $\{h_n: G \to A_n*B_n\}$, be a sequence of homomorphisms.

With the sequence of free products, $\{A_n*B_n\}$, we naturally associate their Bass-Serre trees that we
denote, $\{T_n\}$, with a base point $t_n$ (which is one of the vertices in $T_n$). 
Each homomorphism, $h_n: G \to A_n*B_n$, gives rise to an action $\lambda_{h_n}$
of the group $G$ on
the Bass-Serre tree $T_n$.
 For each index $n$ we fix an element $\gamma_n \in A_n*B_n$, so that the homomorphism $\gamma_n h_n {\gamma_n}^{-1}$
has "minimal displacement", i.e., the element $\gamma_n$ satisfies the equality:
$$\operatornamewithlimits{\max}_{1\le u \le m} \, d_{T_n}(t_n,\gamma_nh_n(g_u)
{\gamma_n}^{-1}(t_n))=
\operatornamewithlimits{\min}_{\gamma \in A_n*B_n}
        \operatornamewithlimits{\max}_{1\le u \le m} \, d_{T_n}(t_n,\gamma h_n(g_u)
{\gamma}^{-1}(t_n))$$
We further set $\mu_n$ to be:
$$\mu_n=\operatornamewithlimits{\max}_{1\le u \le m} \, d_{T_n}(t_n,\gamma_nh_n(g_u)
{\gamma_n}^{-1}).$$

First, suppose that the sequence of integers, $\{\mu_n\}$, is bounded. In that case we can extract 
a subsequence of the homomorphisms $\{h_n\}$ (still denoted $\{h_n\}$), that converges into a limit group (over
free products) $L$, with an associated set of elliptics $E_L$. Furthermore, the sequence of homomorphisms
$\gamma_n h_n \gamma_n^{-1}$ converges into a faithful action of $L$ on some simplicial tree with trivial edge stabilizers,
that we denote $T$. In that case either the entire group $L$ is elliptic (i.e. $E_L=L$), or $L$ is infinite cyclic,
or it is freely decomposable and the stabilizer of each vertex group in $T$ is elliptic. In this case,
the limit group $L$ is a free product of elliptic vertex groups (in $T$) with a (possibly trivial) free group. 

Suppose that the sequence of integers, $\{\mu_n\}$, does not contain a bounded subsequence.
We set ${\{(X_n,x_n)\}}_{n=1}^\infty$ to be the pointed metric spaces obtained 
by rescaling the metric on the Bass-Serre trees $(T_n,t_n)$, by $\mu_n$. 
$(X_n,x_n)$ is endowed with a left isometric action of our f.g.\ group $G$
 via the
homomorphisms ${\gamma_n} h_n \gamma_n^{-1}$.
This sequence of actions of  $G$ on the metric
 spaces 
${\{(X_n,x_n)\}}_{n=1}^\infty$ allows us to obtain an action of $G$ on a real
 tree by passing
to a Gromov-Hausdorff limit.

\vglue 1pc
\proclaim{Proposition 2 ([Pa], 2.3)}  Let $\{X_n\}^\infty
_{n=1}$ be a sequence of $\delta_n$-hyperbolic spaces with
$\delta_\infty \= \varliminf\, \delta_n = 0$.  Let $H$
be a countable group isometrically acting on
$X_n$.  Suppose there exists a base point $x_n$ in $X_n$
such that for every finite subset $P$ of $H$, the sets of geodesics between 
the images of $x_n$ under $P$ form a sequence of 
totally bounded metric spaces.  Then 
there is a subsequence converging in the Gromov topology to a 
$\delta_\infty$-hyperbolic space $X_\infty$ endowed with a 
left isometric action of $H$.
\endproclaim

\smallskip
Our spaces $\{(X_n,x_n)\}^\infty_{n=1}$ endowed with the left isometric
action of $G$, satisfy the 
assumptions of the proposition
  and they are all trees, so they are 0-hyperbolic. Hence,
$X_\infty$ is a real tree endowed with an isometric action of $G$. 
By construction, 
the action of $G$ on the real tree $X_{\infty}$ is non-trivial. Let $\{n_j\}_{j=1}^{\infty}$
be the subsequence for which $\{(X_{n_j},x_{n_j})\}_{j=1}^{\infty}$
 converges to the limit real tree
$X_{\infty}$, and let $(Y,y_0)$ denote this (pointed) limit real tree. 

For convenience, for the rest of this section we (still) denote the homomorphism 
$\gamma_{n_j}h_{n_j}{\gamma_{n_j}}^{-1} : G \to A_{n_j}*B_{n_j}$, by $h_n$.
By passing to a further subsequence we can assume that the sequence of homomorphisms 
$\{h_n\}$ converges into a limit group (over free products) that we denote $L$, with
elliptic elements $E_L$, and an associated quotient map, $\eta:G \to L$, with kernel $K$.
In the sequel we call a limit group over free products that is obtained from a sequence of homomorphisms with
unbounded stretching factors, a $strict$ limit group over free products.

\smallskip
The following simple facts on the kernel of the action, $K$, (see definition 1.1) and the (strict)  limit group $L$
are important observations, and their proof is similar to the proof
of lemma 1.3 of [Se1].

\vglue 1pc
\proclaim{Lemma 3} 
\roster
\item"{(i)}" Elements in $E_L$ fix points in $Y$.
 
\item"{(ii)}" $L$ is f.g.

\item"{(iii)}" If $Y$ is isometric to a real line
 then the limit group $L$ has a subgroup of index at most 2, which is  f.g.\ 
free abelian.

\item"{(iv)}" If $g \in G$ stabilizes a tripod in $Y$ then for all but 
 finitely many $n$'s,  $g \in ker(h_n)$
(recall that a tripod is a finite tree with 3
endpoints). In particular, if $g \in G$ stabilizes a tripod then 
$g \in K$.

\item"{(v)}"  Let $g \in G$ be an element which does not belong to
$K$. Then for all but finitely many $n$'s, $g \notin ker(h_n)$.

\item"{(vi)}"  Every torsion element in $L$ is elliptic, i.e., it is in $E_L$.

\item"{(vii)}" Let $[y_1,y_2] \subset [y_3,y_4]$ be a pair
of non-degenerate segments of $Y$ and assume that the stabilizer of
$[y_3,y_4]$ in $L$,  $stab([y_3,y_4])$, is non-trivial. Then 
$stab([y_3,y_4])$ is an abelian subgroup of $L$ and: 
$$
stab([y_1,y_2])=stab([y_3,y_4])
$$
Hence, the action of $L$ on the real tree $Y$ is (super) stable.

\item"{(viii)}" Let $H<G$ be a f.g.\ subgroup and suppose that $\eta(H) \subset E_L$. Then for
all but finitely many $n$'s, $h_n(H)$ is elliptic, i.e., $h_n(H)$ is contained in a conjugate of
$A_n$ or $B_n$.
\endroster
\endproclaim

\nfp Part (i) follows from the definition of the elliptic elements $E_L$.
A limit group $L$ is a quotient of a f.g.\ group, hence, it is finitely generated. If $Y$ is a real line,
then $L$ contains a subgroup of index at most 2 that acts on the real line by isometries and preserves its
orientation. Hence, this subgroup must be f.g.\ abelian, and it contains no elliptic elements, so it contains
no torsion. Therefore, $L$ contains a f.g.\ abelian subgroup of index at most 2. (iv), (v), and (vii), 
follow by the same argument that is used in the case of free and hyperbolic groups ([Se1],1.3). A torsion
element in $L$ is the image of an element $g \in G$, which is mapped to a torsion element by all the
homomorphisms, $h_n:G \to A_n*B_n$, for large $n$. Hence, $h_n(g)$ is contained in a conjugate of $A_n$ or
$B_n$ for large $n$, so $g$ is mapped to an elliptic element in $L$, and (vi) follows. To prove (viii)
let $H=<u_1,\ldots,u_m>$. Since $H$ is contained in $E_L$ then there exists some $n_0$ so that for all
$n>n_0$, all the elements $u_1,\ldots,u_m$ and the products $u_iu_j$, $i,j=1,\ldots,m$, are mapped to elliptic 
elements by the homomorphism $h_n$. Therefore, for all $n>n_0$, $h_n(H)$ is elliptic, i.e., contained in
a conjugate of $A_n$ or $B_n$. 

\line{\hss$\qed$}

Recall that in limit groups over free and torsion-free hyperbolic groups, every non-trivial abelian subgroup
is contained in a unique maximal abelian subgroup, and every maximal abelian subgroup
is f.g.\ and malnormal. This is clearly not the case in limit groups over free products, as every f.g.\ group
is a limit group over free products. However, for the analysis of strict limit groups over free products,
we are really interested only in non-elliptic abelian subgroups, as only these  
occur as 
stabilizers of non-degenerate segments in real trees on which the strict limit groups act faithfully,
and so that these real trees are obtained as a limit from a sequence of homomorphisms into free products. 
Non-elliptic abelian subgroups have similar properties as abelian subgroups in limit groups over free and torsion-free
hyperbolic groups.
The proof is similar to the proof of lemma 1.4 in [Se1].

\proclaim {Lemma 4} With the notation of definition 1
 let $u_1,u_2,u_3$ be non-trivial elements of 
 $L$, and suppose that at least one of the elements, $u_1,u_2,u_3$, is non-elliptic (i.e., not in $E_L$), 
 and 
$[u_1,u_2]=1$ and $[u_1,u_3]=1$. Then:
\roster
\item"{(i)}"  $u_1,u_2,u_3$ are non-elliptic and $[u_2,u_3]=1$.

\item"{(ii)}" let $A<L$ be  a  non-elliptic abelian subgroup of $L$.
Then $A$ is contained in a unique maximal abelian subgroup in $L$, which is its centralizer, $C(A)$. 
The centralizer of $A$, $C(A)$, intersects the set of elliptic 
elements, $E_L$, trivially.

\item"{(iii)}" let $A$ be a non-elliptic abelian  subgroup in $L$. 
Then the centralizer of $A$, $C(A)$, is almost malnormal in $L$.
$C(A)$ is of index at most 2 in $N(A)$, the normalizer of $A$, and for each element $\ell \in L$, 
$\ell \notin N(A)$, $\ell C(A) \ell^{-1}$  intersects $A$ trivially. Furthermore, if $[N(A):C(A)]=2$ then
$N(A)$ is generated by $C(A)$ and an elliptic element of order 2 that conjugates each element in $C(A)$ to its inverse.
\endroster
\endproclaim

\nfp Let $g_1,g_2,g_3$ be elements of $G$ that are mapped to $u_1,u_2,u_3$. Since $[u_1,u_2]=1$ and
$[u_1,u_3]=1$, and the elements $u_1,u_2,u_3$ are non-trivial, 
there exists some $n_0$, so that for all $n>n_0$, 
$[h_n(g_1),h_n(g_2)]=1$ and $[h_n(g_1),h_n(g_3)]=1$, 
and the elements, $h_n(g_1),h_n(g_2),h_n(g_3)$, are non-trivial. Since for some
$j$, $1 \leq j \leq 3$, $u_j$ is not elliptic, there exists some $n_j>n_0$, so that for all $n>n_j$, $h_n(g_j)$
is a hyperbolic element. Since for $n>n_j$, $h_n(g_1)$ is non-trivial and commutes with $h_n(g_j)$,
$h_n(g_1)$ is a hyperbolic element, and by the same argument so are $h_n(g_2)$ and $h_n(g_3)$. Hence, all
the 3 elements, $h_n(g_1),h_n(g_2),h_n(g_3)$, are hyperbolic and have the same axis, so they all commute and
$[u_2,u_3]=1$. 

By part (i) commutativity is transitive for non-elliptic elements of a limit group over free products. Hence
a non-elliptic abelian subgroup is contained in a unique maximal abelian subgroup, which is its centralizer,
and the centralizer must be non-elliptic as well.

Let $A<L$ be a non-elliptic abelian subgroup. Let $u \in N(A) \setminus C(A)$, and let $g \in G$ be an 
element that is mapped to $u$. Given a finite set of non-trivial elements $g_1,\ldots,g_m$ that are mapped to
$C(A)$, there exists some integer $n_0$, so that for every $n>n_0$, $h_n(g_j)$ are hyperbolic, $h_n(g_i)$ commutes
with $h_n(g_j)$, and $h_n(g)$ does not commute with $h_n(g_i)$,
for all $i,j=1,\ldots,m$, and $h_n(g)h_n(g_j)h_n(g)^{-1}$ commutes with all the elements
$h_n(g_i)$, for $i,j=1,\ldots,m$. This imply that the elements $h_n(g_j)$ have the same axis in the Bass-Serre
tree that is associated with the free product $A_n*B_n$, and the element $h_n(g)$ preserves this axis setwise,
and must be an elliptic element. Hence, $h_n(g)$ is elliptic, and $h_n(g)h_n(g_j)h_n(g)^{-1}=h_n(g_j)^{-1}$, 
$j=1,\ldots,m$. Furthermore, $h_n(g)^2$ is an elliptic element that preserves the axis of the elements
$h_n(g_j)$ pointwise, so $h_n(g)^2=1$. It follows that $ucu^{-1}=c^{-1}$ for every $c \in C(A)$, and $u^2=1$.
By the same argument if $u_1,u_2 \in N(A) \setminus C(A)$ then $u_1u_2 \in C(A)$, hence, $[N(A):C(A)]=2$.

\noindent
Let $\ell \notin N(A)$, and let $t \in G$ be an element that is mapped to $\ell$. Then there exists
some index $n_1$, so that for all $n>n_1$, $h_n(t)$ maps the axis of $h_n(g_1),\ldots,h_n(g_m)$ to a
different axis that intersects the original axis in a bounded (or empty) set. Hence, $\ell C(A) \ell^{-1}$
intersects $C(A)$ trivially.

\line{\hss$\qed$}

Lemma 3 shows that the action of $L$  on the real 
tree $Y$ is (super) stable.
The  original analysis of
 stable actions of  groups on real trees applies to f.p.\ groups
(\cite{Be-Fe1}), and the  limit group $L$ is only
known to be f.g.\, by part (i) of lemma 3.
 Still, given the basic properties
of the action of $L$ on the real tree $Y$ that we already know,
we are able to apply a generalization of Rips' work to f.g.\ groups
obtained in [Se5] and [Gu].  In [Se5] and [Gu],  the real 
tree $Y$ is divided into distinct components, where on each
 component a subgroup of $L$ acts according to one of several canonical types of
actions. The theorem from [Se5] we present, that was later corrected in [Gu],
is going to be used extensively  and its statement uses
the 
 notions and basic definitions that appear in the appendix of
[Ri-Se1]. Hence, we refer a reader who is not yet familiar with these notions
to that appendix and to [Be-Fe1] and [Be].

\vglue 1pc
\proclaim{Theorem 5 (([Se5],3.1),[Gu])}
Let $G$ be a f.g.\ group, let $H_1,\ldots,H_r$ be subgroups of $G$, and suppose that 
$G$ can not be presented as a free product in which the subgroups, $H_1,\ldots,H_r$ can be conjugated 
into the factors. Let $G$ admit  a (super) stable
 isometric action on a real tree 
$Y$, so that  the subgroups, $H_1,\ldots,H_r$, fix points in $Y$. Assume the stabilizer of each 
tripod in $Y$ is trivial. 

\roster
\item"{1)}" There exist canonical orbits of 
subtrees of $Y\colon\  Y_1,\ldots Y_k$ with the following
properties:

\itemitem{(i)}  $gY_i$ intersects $Y_j$ at most in one point if
$i \neq j$.

\itemitem{(ii)}  $gY_i$ is either identical with $Y_i$ or it intersects
it at most in one point.


\itemitem{(iii)}  The action of $\stab (Y_i)$ on $Y_i$ is either discrete or it
 is  of axial 
type or IET type.

\item"{2)}" $G$ is the fundamental group of a graph of groups
 with:

\itemitem{(i)} Vertices corresponding  to orbits of branching points with non-trivial 
stabilizer in
$Y$.  

\itemitem{(ii)} Vertices corresponding to the orbits of the canonical subtrees
$Y_1,\ldots,Y_k$ which are of axial or IET type. The groups associated with these
vertices  are conjugates of the 
stabilizers of these components.  
To a stabilizer of an $IET$ component is an 
 associated 2-orbifold. All boundary components
and branching points in this associated 2-orbifold stabilize points in $Y$. For each
such stabilizer we add edges that connect the vertex stabilized by it and the vertices 
stabilized by its boundary components and branching points.

\itemitem{(iii)} Edges corresponding to orbits of edges between branching points 
with non-trivial stabilizer in
the discrete part of $Y$ with
edge groups which are conjugates of the stabilizers
of these edges.

\itemitem{(iv)} Edges corresponding to orbits of points of intersection
between the orbits of $Y_1,\ldots,Y_k$. 
\endroster
\endproclaim

\smallskip
Before concluding our preliminary study of limit groups over free products that appear as
limits of sequences of homomorphisms with unbounded stretching factors $\{\mu_n\}$, and their actions
on the limit real tree, we present the following basic facts that are 
necessary in the sequel.

\proclaim {Proposition 6} Suppose that  $L$ is a (strict) limit group over free products, 
that is obtained as a limit of
homomorphisms into free products with unbounded stretching factors $\{\mu_n\}$. $E_L$ is its set
of elliptic elements,  and the limit real tree on
which $L$ acts that is obtained from this sequence of homomorphisms is $(Y,y_0)$. Suppose further that $L$ 
does not admit a non-trivial free decomposition in which all the elements in the set $E_L$  can
be conjugated into the various factors. Then:
\roster
\item"{(i)}" Stabilizers of non-degenerate segments which lie in the 
complement of the discrete and axial parts of  $Y$ are trivial in $L$.

\item"{(ii)}"  
The (set) stabilizer of an axial component in $Y$ is
either a maximal abelian subgroup in $L$, or it contains a maximal abelian subgroup in $L$ as a 
subgroup of index 2.

\item"{(iii)}" Let $A$ be the maximal abelian subgroup that is contained in the set stabilizer of an
axial component in $Y$. $A$ 
can be presented as a direct sum $A=A_1+A_2$, 
where $A_1$ is the (possibly trivial) point stabilizer of the axial component, and $A_2$ is
a f.g\ free abelian group that acts freely
on the axial component, and $A_2$ has rank at least 2.
\endroster
\endproclaim

\nfp Since the elements in $E_L$ fix points in the limit tree $Y$ (Part (i) in lemma 3), 
the action of $L$ on the real
tree $Y$ satisfies the conclusion of theorem 5. Hence, the stabilizer of a segment in an IET component
in $Y$ fixes the entire IET component, and in particular it fixes a tripod. By  part (iv) of lemma 3
a stabilizer of a tripod is trivial, hence, so is the stabilizer of a non-degenerate segment in an
  IET component in $Y$. 

Let $Ax$ be an axis of an axial component in $Y$, and let $stab(Ax)$ be its set stabilizer. Let $stab^{+}(Ax)$
be the subgroup of $stab(Ax)$ that preserve the orientation of $Ax$. By the same argument that
was used in the proof of lemma 4, $stab^{+}(Ax)$, is abelian. Since $stab(Ax)$ normalizes $stab^{+}(Ax)$,
lemma 4 implies that the index of $stab^{+}(Ax)$ in $stab(Ax)$ is bounded by 2.

Let $A=stab^{+}(Ax)$, and let $A_1<A$ be the point stabilizer of $Ax$. Then, by theorem 5 (the structure of an axial
component) there exists a short exact sequence:
$1 \to A_1 \to A \to B \to 1$, where $B$ is a f.g.\ free abelian group. Since $A$ is abelian and $B$ is free 
abelian, the short exact sequence splits, and $A$=$A_1$+$A_2$, where  $A_2$ is
isomorphic to $B$, hence, $A_2$ is f.g.\ free abelian.

\line{\hss$\qed$}

By theorem 5  and proposition 6 
a non-trivial strict limit group over free products, which is not a  cyclic group, 
admits a non-trivial virtually abelian
decomposition (i.e., a graph of groups with virtually abelian edge groups).
 To further study the algebraic structure of
a strict limit group we need to construct its canonical virtually 
abelian JSJ decomposition. However, unlike the case of
limit groups over free and hyperbolic groups, in constructing the virtually abelian
JSJ decomposition of a strict limit group over free products, we will not be interested in all the virtually 
abelian 
decompositions of $L$, but only in those virtually
abelian decompositions in which all the elements in $E_L$ are elliptic (i.e., can be conjugated into
non-virtually-abelian, non-QH vertex groups), and for which the (non-trivial) maximal 
abelian subgroups that are contained in 
the virtually abelian
edge groups are not in $E_L$.
Note that since a non-trivial strict limit group over free products admits a virtually abelian decomposition 
in which the elements $E_L$ can be conjugated into non-QH, non-virtually-abelian vertex groups, 
and the (non-trivial) maximal abelian subgroups of edge groups are not in $E_L$,
if a strict limit group over free products is not virtually abelian nor a Fuchsian group, 
its (virtually) abelian JSJ decomposition is non-trivial.

To construct the virtually abelian JSJ decomposition of a strict limit group over free products we need to study
some basic properties of virtually abelian splittings. We start with the following lemma,
which is identical to lemma 2.1 in [Se1] (the proofs are identical as well).

\vglue 1pc
\proclaim{Lemma 7} Let $L$ be a  strict limit group over free products with set of elliptics $E_L$, and suppose
that $L$ admits no free product in which the elements in $E_L$ can be conjugated into
the various factors. Let $A$ be a  non-elliptic abelian subgroup in $L$, 
and let $M$ be the normalizer of $A$ in $L$. Suppose that $M$ is abelian. Then: 
\roster
\item"(i)" If $L=U*_A V$, and the elements in $E_L$ can be conjugated into $U$ and $V$,  and $M$ is not cyclic,
then $M$ can be conjugated into either $U$ or $V$.

\item"(ii)" If $L=U*_A$, and the elements in $E_L$ can be conjugated into $U$,  and $M$ is not cyclic,
 then either $M$ can be conjugated into $U$, or $M$ can be conjugated
 onto $M'$, so that  $L=U*_A{M'}$. 
\endroster
\endproclaim

Unlike limit groups over free and torsion-free hyperbolic groups in which normalizers of non-trivial
abelian subgroups are abelian, by proposition 6 
the normalizers of non-elliptic abelian subgroups in $L$
are either abelian or virtually abelian, and if they are not abelian, the abelian 
centralizers of these (non-elliptic) abelian subgroups 
are contained in their normalizers as  subgroups of index 2. Lemma 7 deals with the case in which the normalizer of
such an abelian subgroup is abelian. To construct the JSJ decomposition of limit groups over free products, 
we still need
to analyze splittings over non-elliptic abelian subgroups with virtually abelian, non-abelian normalizers. 

\vglue 1pc
\proclaim{Lemma 8} Let $L$ be a limit group over free products,
and let $A$ be a non-elliptic abelian subgroup in $L$. Let $E_L$ be the set of elliptics in $L$,
and suppose that $L$ admits no free product decomposition in which the elements of $E_L$ can be conjugated into
the factors. Let $C(A)$ be the centralizer
of $A$, let $M$
be the normalizer of $A$, and suppose that $[M:C(A)]=2$. Then: 
\roster
\item"(i)" If $L=U*_A V$, and all the elements in $E_L$ can be conjugated into $U$ or $V$, 
then either $M$ can be conjugated into either $U$ or $V$, or $M$ can be conjugated onto $M'$, and $M'$
inherits an abelian 
amalgamation: $M'=U_1*_{A}V_1$, where $U_1<U$, $V_1<V$, $[U_1:A]=[V_1:A]=2$, and both $U_1$
and $V_1$ are generated by $A$ and an element of order 2. In this case, $M$ is the semidirect
product of $A$ with an infinite dihedral group. In this case we can modify the given abelian decomposition,
and obtain a virtually abelian decomposition, $L=U*_{U_1} \, M'*_{V_1} \, V$.

\item"(ii)" If $L=U*_A$, and the elements in $E_L$ can be conjugated into $U$, 
then either $M$ can be conjugated into $U$, or $M$ can be conjugated onto $M'$, and $M'$
inherits an abelian 
amalgamation: $M'=U_1*_{A}V_1$, where $U_1<U$, $V_1<tUt^{-1}$, where $t$ is a Bass-Serre generator that
is associated with the splitting, $L=U*_A$. $[U_1:A]=[V_1:A]=2$, and both $U_1$
and $V_1$ are generated by $A$ and an element of order 2. In this case, $M$ is the semidirect
product of $A$ with an infinite dihedral group. In the HNN case, $L=U*_A$, 
we can modify the given abelian decomposition,
and obtain a virtually abelian decomposition, $L=(U*_{U_1} \, M')*_{V_1}$, where with $V_1$ there are two associated 
embeddings, one into $M'$ and one into $tUt^{-1}$. The graph of groups that is associated with this decomposition
contains two vertices (with vertex groups, $U$ and $M'$), and two edges with edge groups, $U_1$ and $V_1$.
\endroster
\endproclaim

\nfp Let $L=U*_A V$ and suppose that $M$, the normalizer of $A$, that contains the centralizer of $A$ as a 
subgroup of index 2,   is not elliptic. Then $M$ preserves (setwise) an axis in the Bass-Serre tree that is
associated with the amalgamated product $U*_A V$. Since $A$ preserves this axis pointwise, and $M$ contains an
(elliptic) element that acts on the axis as an inversion, $M/A$ acts on the axis as an infinite  dihedral group.
Hence, it inherits from it a splitting: $M=U_1*_A V_1$, where $U_1$ and $V_1$ contain $A$ as a subgroup of index 2,
and they are both obtained from $A$ by adding to it an elliptic element of order 2. If we start with
the Bass-Serre tree that is associated with $U*_A V$, add a vertex
in the middle of the edge that is stabilized by $A$ and connected the vertices that are stabilized by $U$ and $V$,
and then fold the couple of edges that are stabilized by $A$ and associated with the elements of order 2
in $U_1$ and $V_1$, we obtain the splitting:
$L=U*_{U_1} \, M'*_{V_1} \, V$, i.e., a splitting of $L$ with two vertex groups $U$ and $M'$, and two
edge groups, $U_1$ and $V_1$. 

Let $L=U*_A V$. the argument that we use in this case is similar.
Suppose that $M$, the normalizer of $A$,  
is not elliptic. In this case, $M$ preserves (setwise) an axis in the Bass-Serre tree that is
associated with the HNN extension $U*_A$. $A$ preserves this axis pointwise, and $M$ contains an
(elliptic) element that acts on the axis as an inversion, hence, $M/A$ acts on the axis as an infinite  dihedral group.
Therefore, as in the amalgamated product case, $M$ inherits from this action a splitting: $M=U_1*_A V_1$, 
where $U_1$ and $V_1$ contain $A$ as a subgroup of index 2,
and they are both obtained from $A$ by adding to it an elliptic element of order 2. $U_1<U$, and $V_1<tUt^{-1}$, for
an appropriate Bass-Serre generator $t$.  If we start with
the Bass-Serre tree that is associated with $U*_A $, add a vertex
in the middle of the edge that is stabilized by $A$ and connects the vertices that are stabilized by $U$ and $tUt^{-1}$,
and then fold the couple of edges that are stabilized by $A$ and associated with the elements of order 2
in $U_1$ and $V_1$, we obtain the splitting:
$L=(U*_{U_1} \, M')*_{V_1}$, where $V_1$ embeds into $V_1$ and into $tUt^{-1}$, i.e., the limit group $L$ admits
a graph of groups decomposition with two vertex groups, stabilized by $U$ and $M'$, and two edges in between
these two vertices, one edge is stabilized by $U_1$ and the second by $V_1$.

\line{\hss$\qed$}

According to lemma 7  we replace each abelian splitting of $L$ of
the form $L=U*_A$, in which all the elements in $E_L$ can be conjugated into $U$, 
$A$ is a non-elliptic abelian subgroup 
in $L$, and the centralizer of $A$ which is also its normalizer is denoted by $M$, and $M$ can not be conjugated 
into $U$,
 by the amalgamated product 
$L=U*_A M'$ (where $M'$ is a conjugate of $M$). 
According to part (i) of lemma 8  we replace each abelian splitting of $L$ of
the form $L=U*_A V$, in which all the elements in $E_L$ can be conjugated into $U$ and $V$,
$A$ is a  non-elliptic subgroup in $L$, and $M$ the normalizer of $A$ contains the centralizer of $A$ as
a subgroup of index 2, and $M$ can not be conjugated into $U$ nor $V$,
 by the amalgamated product 
$L=U*_{U_1} \, M'*_{V_1} \, V$, where $M'$ is a conjugate of $M$. If $L=U*_A$, the elements in $E_L$ can be 
conjugated into $U$, $A$ is non-elliptic in $L$,  the normalizer $M$ of $A$ in $L$ contains the centralizer
of $A$ as a subgroup of index 2, and $M$ can not be conjugated into $U$, then we replace the given HNN extension
by a graph of groups with two vertices and two edges between them, according to part (ii) of lemma 8, 
$L=(U*_{U_1} \, M')*_{V_1}$. 

By performing these replacements, we get that every non-elliptic abelian subgroup of
$L$ with a non-cyclic centralizer is contained in a vertex group
 in all the virtually abelian splittings of $L$ under consideration,
i.e., splittings in which edge groups are non-elliptic abelian,  or edge groups that 
contain non-elliptic abelian subgroups
as subgroups of index 2,
and all the elements in $E_L$ can be conjugated into
the  vertex groups.
This will allow us to use acylindrical accessibility in analyzing all the
 abelian splittings of the limit group over free products $L$ that are obtained from converging sequences
of homomorphisms into free products.

\vglue 1pc
\proclaim{Definition 9 ([Se5],[De],[We])} A splitting of a group $H$ is called 
$k$-$acylindrical$ if for every element $h \in H$ which is not the
identity, the fixed set of $h$ when acting on the Bass-Serre tree
corresponding to the splitting has diameter at most $k$. Following Delzant [De],
 and Weidmann [We], we say that a splitting of $H$ is $(k,C)$-acylindrical if the stabilizer of a
path of length bigger than $k$ in the Bass-Serre tree corresponding to the splitting has stabilizer of
cardinality at most $C$.
\endproclaim

If a strict limit group over free products $L$ can be written in the form $L=V_1*_{A_1} \, V_2*_{A_2} \, V_3*_{A_3} \, V_4$, where 
$A_1,A_2,A_3$ 
 are subgroups of a maximal abelian subgroup $M$, which is its own normalizer, and $M$ is a subgroup of $V_1$, 
then one can modify the 
corresponding graph of groups to a tripod of groups with $V_1$ in the center, and $V_2,V_3,V_4$
at the 3 roots. Similarly if $A_1,A_2,A_3$ are subgroups of a maximal abelian subgroup, which is of index 2 in its
normalizer $M$, then if $M$ is contained in 
one of the vertex groups $V_i$, $i=1,\ldots,4$, then
one can modify the corresponding graph of groups to a tripod of groups in the same way. If $M$ is not contained in
one of the limit groups, $V_i$, $i=1,\ldots,4$, then one can modify the given splitting to a virtually abelian splitting
which is a tree with one root and 4 vertices connected to it, where $M$ is placed at the root, and the subgroups
$V_i$, $i=1,\ldots,4$, are placed at the vertices that are adjacent to the root.

Since by lemma 4 the centralizer of a non-elliptic abelian subgroup
 is almost malnormal,
the Bass-Serre trees corresponding to these tripods and quadruple of groups are (2,2)-acylindrical.
This folding and sliding operation generalizes to an arbitrary (finite) virtually abelian splitting 
of a limit group over 
free products over normalizers of non-elliptic abelian subgroups, 
where all the elements in $E_L$ can be conjugated into
non-abelian, non-QH vertex groups.

\vglue 1pc
\proclaim{Lemma 10} Let $L$ be a limit group over free products that does not admit a free splitting in which
all the elements in $E_L$ can be conjugated into the various factors.
A splitting of $L$, in which all the edge groups 
are non-elliptic abelian subgroups in $L$, and in which all the elements in $E_L$ can be conjugated into non-QH, 
non-abelian vertex groups, 
can be modified using lemmas 7 and 8 to a virtually abelian splitting (of $L$)
in which all the normalizers of non-elliptic abelian subgroups with non-cyclic centralizers 
can be conjugated into non-QH vertex groups,
 and so that the obtained virtually
abelian splitting is $(2,2)$-acylindrical. 
\endproclaim

\nfp Let $L$ be a limit group over free products that admits no free factorization in which the
elements of $E_L$ can be conjugated into the factors. Let $\Lambda$ be a graph of groups with fundamental
group $L$ with non-elliptic abelian edge groups. 
If the normalizer of an abelian edge group in $\Lambda$ can not
be conjugated into a vertex group in $\Lambda$, we perform the modification that appears in part (ii)
of lemma 7 in case the normalizer of an edge group is abelian, and the modification of parts (i) and (ii)
in lemma 8 in case
the normalizer of an edge group is not abelian, where these modifications are applicable. After performing 
these modifications, and sliding operations, so that the fixed set of a non-elliptic abelian subgroup 
is star-like, we obtain a graph of groups $\Lambda'$, with virtually abelian
edge groups, in which every non-elliptic abelian subgroup with non-cyclic centralizer can be conjugated into
a non-QH vertex group. 
In the Bass-Serre
tree that corresponds to $\Lambda'$ the fixed set of every non-elliptic element has diameter bounded by 2.
Since the centralizers of non-elliptic abelian subgroups are almost malnormal, the stabilizers of paths of
length 3 in $\Lambda'$ are either trivial, or a cyclic subgroup of order 2. Hence,  $\Lambda'$ is
(2,2)-acylindrical.

\line{\hss$\qed$}

Lemma 10 shows that if in all the virtually abelian splittings of $L$ 
under consideration, 
all the normalizers of non-cyclic abelian subgroups can be conjugated into vertex groups, 
these virtually abelian splittings are 
$(2,2)$-acylindrical.  
 This acylindricity finally enables one to construct the canonical 
virtually abelian JSJ decomposition of a strict limit group over free products (see section 2 of [Se1]).

\vglue 1pc
\proclaim{Theorem 11 (cf. ([Se1],2.7))} Suppose
 that $L$ is a strict limit group over free products with set of elliptics $E_L$,
so that $L$ admits no free decomposition in which the elements in $E_L$ can be conjugated into the various
factors.
 There exists a  reduced unfolded
splitting of $L$ with virtually abelian edge groups, 
which we  call the virtually abelian   $JSJ$ (Jaco-Shalen-Johannson)
$decomposition$
 of $L$, with the
 following properties:

\roster
\item"{(i)}" Every  canonical
  maximal QH subgroup (CMQ) of $L$ is
conjugate to a vertex group in the  JSJ decomposition. Every QH
subgroup of $L$, in which all the elements in $E_L$ can be conjugated into vertex groups that are adjacent to the
$QH$ subgroup or into torsion elements in the QH subgroup, can be conjugated into one of the CMQ subgroups of $L$.
Every vertex group in the JSJ decomposition which is not 
 a CMQ subgroup of $L$ is elliptic in any
abelian splitting of $L$ under consideration. 

\item"{(ii)}" Every CMQ subgroup is a Fuchsian group (in general, with torsion), where all its torsion 
elements are elliptic in $L$. The edge groups that are connected to a CMQ subgroup, that  are all  cyclic,
may be  elliptic. 

\item"{(iii)}" Every edge group that is not connected to a CMQ vertex group in the JSJ decomposition, or an
edge group that is connected to a virtually abelian vertex group, contains
an abelian subgroup of index at most 2, and this abelian subgroup is non-elliptic. 

\item"{(iv)}" A one edge abelian splitting $L=D*_AE$ or
 $L=D*_A$, in which $A$ is a non-elliptic abelian subgroup, and all the elements in $E_L$ can be conjugated into
$D$ or $E$, 
which  is hyperbolic in another  such elementary
abelian splitting,
is obtained from the virtually abelian JSJ decomposition of $L$  by cutting
 a 2-orbifold corresponding
 to a CMQ
subgroup of $L$ along a  weakly essential s.c.c.\.

\item"{(v)}" Let $\Theta$ be a one edge splitting of $L$ along a  non-elliptic abelian
subgroup,
 $L=D*_AE$ or
$L=D*_A$, 
in which all the elements in $E_L$ can be conjugated into $D$ or $E$.
  Suppose that the given elementary splitting is elliptic with respect to any other such
elementary abelian splitting of $L$.
Then $\Theta$ is obtained from the
JSJ decomposition of $L$ by a sequence of collapsings, foldings,
conjugations, and finally possibly unfoldings that reverse the foldings that are 
performed according to part (i) of lemma 7 and parts (i) and (ii) of lemma 8. 


\item"{(vi)}" If $JSJ_1$ is another JSJ decomposition of $L$,
 then $JSJ_1$
is obtained from the JSJ decomposition by a sequence of slidings, conjugations
and modifying boundary monomorphisms by conjugations (see section 1 of
[Ri-Se2] for these notions)
\endroster
\endproclaim

\nfp By lemma 10 the splittings of the ambient limit group (over free products) $L$ that are considered 
for the construction of the virtually abelian JSJ decomposition of $L$, have the property that all the elliptic elements
$E_L$ in $L$ can be conjugated into non-QH, non-abelian vertex groups, and every non-cyclic, non-elliptic
abelian subgroup of $L$ can also be conjugated into a non-QH vertex group. Since $L$ admits no free decompositions in
which the elements $E_L$ can be conjugated into the factors, there is no pair of elliptic-hyperbolic splittings
along non-elliptic abelian groups (so that the elements $E_L$ can be conjugated into vertex groups), i.e., 
all the splittings along non-cyclic, non-elliptic abelian groups under consideration,  
are elliptic-elliptic with respect to all the splittings along  non-elliptic abelian groups under consideration.

Since the modifications of abelian splittings along non-elliptic abelian subgroups that are performed according
to lemmas 7 and 8, are performed only in case the centralizers of (non-elliptic) edge groups are non-cyclic,
in case centralizers of non-elliptic edge groups are (infinite) cyclic, we consider only cyclic edge groups
(and not dihedral ones). Hence, the only hyperbolic-hyperbolic splittings under consideration are 
pairs of splittings along infinite cyclic groups. For these we can apply [Ri-Se2] (this part of [Ri-Se2],
the construction of the quadratic part (section 5 in the paper),
applies to f.g.\ groups, and not only to f.p.\ ones), that produces a finite collection of CMQ subgroups of $L$,
and a quadratic decomposition of $L$ with the CMQ subgroups as part of the vertex groups, 
so that every splitting of $L$ along a non-elliptic cyclic group, in which
the elements of $E_L$ can be conjugated into vertex groups, and so that this splitting is hyperbolic
with respect to another such splitting, is obtained from the quadratic decomposition of $L$ by cutting one of
the CMQ subgroups along a s.c.c.\ and possibly collapsing the rest of the splitting.

Given the quadratic decomposition of $L$, to complete the construction of the (virtually abelian) JSJ decomposition
of $L$, we successively refine the quadratic decomposition using splittings that are elliptic with respect
to it (and in which the elements of $E_L$ can be conjugated into vertex groups). This refinement process
terminates after finitely many steps, since all the obtained splittings are (2,2) acylindrical, and by
[We] this implies a bound on the combinatorial complexity of the obtained splitting. All the properties of the
obtained (virtually) abelian JSJ decomposition of $L$ follow in the same way as in section 7 of [Ri-se]. 

\line{\hss$\qed$}

\vglue 1.5pc
\centerline{\bf{\S2. A Descending Chain Condition}}
\medskip

In section 4 of [Se1] we were able to use the cyclic JSJ decomposition of a ($F_k$) limit group,
in order to show that ($F_k$) limit groups are f.p.\ and that a f.g.\ group is a limit group
if and only if it is $\omega$-residually free. For limit groups over a torsion-free
hyperbolic group, we were able to prove similar d.c.c.\ and a.c.c.\ as in the case of a free group,
even though a limit group over hyperbolic groups need not be finitely presented.

Limit groups over free products do not satisfy the d.c.c.\ and a.c.c.\ conditions that limit groups over
free and hyperbolic groups do satisfy. However, weaker principles do hold for these limit groups, and these
are enough to construct Makanin-Razborov diagrams, that encode sets of solutions to systems of equations over
free products. As we will see one of the keys to formulate and prove the d.c.c.\ and a.c.c.\ principles
that we  present for limit groups over free products, is our consideration of limit groups over the entire
class of free products, and not over a given one.

We start with a d.c.c.\ for limit groups over free products which is a key to our 
entire approach. It uses the techniques that were used to prove a general d.c.c.\ for limit
groups over hyperbolic groups, but it is not as general as in the case of limit groups over free
and hyperbolic groups. 

\vglue 1pc
\proclaim{Definition  12} Let $G$ be a f.g.\ group.
On the set of  limit groups over free products that are quotients of $G$, together with the quotient maps
from $G$ to these limit groups, we define a partial order.
Let $L_1, L_2$  be two limit groups over free products that are quotients of $G$,
with sets of elliptics, $E_{L_1}, E_{L_2}$, in correspondence,
and with prescribed quotient  maps $\eta_i: G \to L_i$, $i=1,2$. We write that  $(L_1,\eta_1) > (L_2,\eta_2)$,
if there exists an  epimorphism: $\tau: L_1 \to L_2$, that maps the elliptics in 
$L_1$ into the elliptics of $L_2$, $\tau(E_{L_1}) \subset E_{L_2}$,
and for which $\tau:L_1 \to L_2$ satisfies either:
\roster
\item"{(1)}" $\tau$  has a non-trivial kernel.

\item"{(2)}" $\tau$ is an isomorphism, and $\tau(E_{L_1})$ is a proper subset of $E_{L_2}$.  
\endroster

If there exists an isomorphism $\tau:L_1 \to L_2$ that maps the elliptics in $L_1$
onto the elliptics in $L_2$, and for which: $\eta_2=\tau \circ \eta_1$, 
we say that $(L_1,\eta_1)$ is in the same equivalence class as $(L_2,\eta_2)$. Note that the relation that is 
defined on the limit quotients (over free products) of a f.g.\ group is a partial order. 
\endproclaim

\vglue 1pc
\proclaim{Theorem 13} Let 
$G$ a f.g.\ group. Every strictly decreasing sequence of limit groups over free products 
that are quotients of $G$:
$$(L_1,\eta_1) \, > \, (L_2,\eta_2) \, > \, (L_3,\eta_3) \, > \, \ldots$$
for which:
\roster
\item"{(1)}" the maps: $\tau_i:L_{i} \to L_{i+1}$, that satisfy:
$\eta_{i+1}=\tau_i \circ \eta_i$, are proper quotient maps (i.e., have non-trivial kernels).

\item"{(2)}" the maps $\tau_i$ do not map non-trivial elements in $E_{L_i}$ to the identity element in $L_{i+1}$.
\endroster
terminates after finitely many steps.
\endproclaim

\nfp The argument that we use is a modification of the argument that is used to prove theorem 1.12 in [Se3].
Suppose that there exists a f.g.\ group $G$ 
for which there exists an infinite decreasing sequence of
limit groups over free products that are quotients of $G$: 
$L_1>L_2>L_3> \ldots$ that satisfy the conditions of the theorem. 
W.l.o.g. we may assume that the f.g.\ group $G$ is a free group $F_d$, for some integer $d$.
We fix $F_d$, where $d$ is the minimal positive integer for which there exists
an infinite descending chain of limit groups over free products so that consecutive quotient maps, 
$\tau_i:L_i \to L_{i+1}$,
have non-trivial kernels and do not map non-trivial elliptic elements to the identity element, 
and  fix a free basis for $F_d$,  
$F_d=<f_1,\ldots,f_d>$. We set $C$ to be the Cayley graph of
$F_d$ with respect to the given generating set, and look at an infinite 
decreasing sequence constructed in the following way. We set $R_1$ to be a limit group over
free products, which is a quotient of $F_d$,
with the following properties:
\roster
\item"{(1)}" $R_1$ is a proper quotient of $F_d$.

\item"{(2)}" $R_1$ can be extended to an infinite decreasing sequence of
 limit groups over free groups: $R_1 > L_2>L_3> \ldots$, that satisfy the conditions of the theorem.

\item"{(3)}" The map $\eta_1: F_d \to R_1$ maps to the identity the maximal number of elements in the ball of radius 1 in the Cayley graph $C$,
 among all possible maps 
from $F_d$ to a  limit group over free products $L$, that satisfy properties (1) and (2).
\endroster
We continue iteratively. At step $n$, given the finite decreasing sequence
$R_1>R_2 > \ldots > R_{n-1}$, we choose the  limit group over free products,
$R_{n}$, to satisfy:
\roster
\item"{(1)}" $R_n$ is a proper quotient of $R_{n-1}$.

\item"{(2)}" The finite decreasing sequence of  limit groups over free products:
$R_1>R_2> \ldots > R_n$ can be
extended to an infinite decreasing sequence that satisfies the conditions of the theorem.

\item"{(3)}" The map $\eta_n: F_d \to R_n$ (that is obtained as a composition of the map $F_d \to R_1$
with the sequence of proper epimorphisms:
$R_i \to R_{i+1}$, $i=1,\ldots,n-1$, 
maps to the identity the maximal number of elements in the ball of radius $n$ in the Cayley graph $C$,
 among all the possible maps 
from $F_d$ to a  limit group over free products, $L_n$, that satisfy 
the properties (1) and (2).
\endroster

To prove theorem 13, we will show that the last descending sequence we constructed terminates after
finitely many steps. 
With the decreasing sequence $R_1>R_2> \ldots$ we associate a sequence of
homomorphisms into free products: $\{h_n: F_d \to A_n*B_n \}$.
For each index $n$, $R_n$ is a quotient of $F_d$, 
 hence, $R_n$ is generated by $d$ elements that are the image of the fixed generators of
$F_d$ under the quotient map $\eta_n$. 

\noindent
$R_n$ is a  limit group over free products. Hence,  $R_n$, with its set of elliptics $E_{R_n}$,
is obtained from a convergent sequence of homomorphisms
$\{ u_s : G_n \to C_s*D_s \}$, where $G_n$ is a f.g.\ group. 
Since $R_n$ is generated by the image of the elements $f_1,\ldots,f_d$
under the quotient map $\eta_n$, for large enough $s$, the images $u_s(G_n)$ are d-generated
groups, and furthermore, they are generated by the images of $d$ elements in 
the f.g.\ group $G_n$, that are mapped by the quotient map $\nu_n: G_n \to R_n$ onto the elements
$\eta_n(f_1),\ldots,\eta_n(f_d)$. Hence, we may assume that the limit groups over free products,
$R_n$, are obtained as limit groups from a sequence of homomorphisms 
$\{v_s:F_d \to C_s*D_s \}$, and the image of the fixed generating set of the free group $F_d$,
is the set of elements  
$\eta_n(f_1),\ldots,\eta_n(f_d)$. 

For each index $n$,
we pick $h_n$ to be a homomorphism $h_n: F_d \to A_n*B_n$, so that $h_n$ is a homomorphism
$v_s:F_d \to C_s*D_s$ for some large index $s$, so that $h_n$ satisfies the following two
conditions: 
\roster
\item"{(i)}" every element in the ball of radius $n$ of $C$, the Cayley graph of $F_d$, that is mapped 
 by the quotient map $\eta_n: F_d \to R_n$ to the trivial 
element, is mapped by $h_n$ to the trivial element in $A_n*B_n$. Every
such element that is mapped to a non-trivial element by $\eta_n$, is mapped
by $h_n$ to a non-trivial element in $A_n*B_n$.

\item"{(ii)}" every element in the ball of radius $n$ of $C$, the Cayley graph of $F_d$, that is mapped 
 by the quotient map $\eta_n: F_d \to R_n$ to an elliptic element, 
is mapped by $h_n$ to an elliptic element in $A_n*B_n$. Every
such element that is mapped to a non-elliptic element by $\eta_n$, is mapped
by $h_n$ to a non-elliptic element in $A_n*B_n$.
%
\endroster

From the sequence $\{h_n\}$ we can extract a subsequence that converges
into a limit group over free products, that we denote $R_{\infty}$. By construction, the limit
group $R_{\infty}$ is the direct limit of the sequence of (proper) epimorphisms:
$$F_d \to R_1 \to R_2 \to \ldots$$

Let $\eta_{\infty}: F_d \to R_{\infty}$ be the canonical quotient map.
Our approach towards  proving the termination of given descending chains of limit groups over free products
is based on studying the structure of the limit group $R_{\infty}$, and its
associated quotient map $\eta_{\infty}$. We
start this study by listing some basic properties of them.

\vglue 1pc
\proclaim{Lemma 14} 
\roster
\item"{(i)}" $R_{\infty}$ is not finitely presented. 

\item"{(ii)}" $R_{\infty}$ can not be presented as the free product of a f.p.\ group and 
freely indecomposable elliptic subgroups. 

\item"{(iii)}" Let
$R_{\infty}=U_1*\ldots*U_t*F$ be the most refined (Grushko) free decomposition of $R_{\infty}$ in which the
elliptic elements in $R_{\infty}$, $E_{R_{\infty}}$, can be conjugated into the various factors, and
 $F$ is a f.g.\ free group. Then
there exists an index $j$, $1 \leq j \leq t$, for which:
\itemitem{(1)} $U_j$ is not finitely presented nor elliptic.

\itemitem{(2)} If $B$ is a f.g.\ subgroup of $F_d$ for which $\eta_{\infty}(B)=U_j$, then $U_j$
is a strict limit group over free products of a subsequence of the restricted homomorphisms,
$h_n|_B$. 
\endroster
\endproclaim

\nfp To prove part (i), suppose that $R_{\infty}$ is f.p. i.e.: 
$$R_{\infty} \ = \ <g_1,\ldots,g_d \, | \, r_1,\ldots,r_s>.$$
Then for some index $n_0$, and every index $n>n_0$, $h_n(r_j)=1$ for $j=1,\ldots,s$. This implies
that for some index $n_1>n_0$, and every index $n>n_1$, each of the  groups $R_n$ is a
quotient of $R_{\infty}$, by a quotient map that send the generating set $g_1,\ldots,g_d$ of
$R_{\infty}$ to the elements $\eta_n(f_1),\ldots,\eta_n(f_d)$, a contradiction.

Suppose that $R_{\infty}=V_1* \ldots *V_t *M$ where $M$ is f.p. and each of the factors $V_j$ is 
elliptic.
Let $B_1,\ldots,B_t$ and $D$ be f.g.\
subgroups of $F_d$, for which $\eta_{\infty}(B_j)=V_j$ for $j=1,\ldots,t$, and $\eta_{\infty}(D)=M$.
W.l.o.g. we may assume that the free group $F_d$ is generated by the collection of the subgroups
$B_1,\ldots,B_t,D$.

\noindent
Since the factors $V_j$, $j=1,\ldots,t$, are elliptic, and since the subgroups, $B_j$, $j=1,\ldots,t$,
are f.g.\ 
 for every index $j$, $j=1,\ldots,t$, there exists an index 
$n_j$, so that for every index $n>n_j$, the image $\eta_n(B_j)$ is elliptic. Since the maps $\tau_i:R_i \to R_{i+1}$
do not map non-trivial  elliptic elements (in $E_{R_i}$) to the identity element,  $\eta_n(B_j)$ is isomorphic to 
$\eta_{\infty}(B_j)=V_j$ via the map $\eta_{\infty} \circ \eta_n^{-1}$. 

\noindent
The factor $M$ is assumed f.p., hence, if $D=<d_1,\ldots,d_s>$, then $M=<d_1,\ldots,d_s \, | \,
r_1,\ldots,r_u>$. There exists an index $n_0$, for which for every index $n>n_0$, $\eta_n(r_i)=1$,
for $i=1,\ldots,u$.

Let $m_0>n_j$ for $j=0,\ldots,t$. By our arguments, from the universality of free products,
all the relations that hold in $R_{\infty}$ hold in $\eta_{m_0}(F_d)=R_{m_0}$. Hence, $R_{m_0}$ is a
quotient of $R_{\infty}$, where the quotient map maps the prescribed generators of $R_{\infty}$ to the
prescribed generators of $R_{m_0}$ (i.e., the corresponding images of the given set of generators
$F_d=<f_1,\ldots,f_d>$). Since $R_{m_0+1}$ is a proper quotient of $R_{m_0}$, this implies that $R_{m_0+1}$
is a proper quotient of $R_{\infty}$, again by a map that maps the prescribed set of generators of $R_{\infty}$
to the prescribed set of generators of $R_{m_0+1}$, which clearly contradicts our assumptions that the sequence
of limit groups $\{R_j\}$ is strictly decreasing with $R_{i+1}$ being a proper quotient of $R_i$ for every
index $i$,  and the limit group $R_{\infty}$ is the direct limit of this decreasing sequence. This concludes the proof of 
part (ii). 

To prove part (iii) note that (1) in part (iii) follows from part (ii). Every factor $U_j$ of the limit group $L$
that is not elliptic is a strict limit group that is obtained from a sequence of homomorphisms of some
f.g\ subgroup of $F_d$, and (2) follows.

\line{\hss$\qed$}

$R_{\infty}$ is a limit group over free products which  is a (proper) quotient of all the limit groups over 
free products,
$\{R_n\}$. 
For each index $n$, the limit group $R_n$ was  chosen to maximize the number of elements 
that are mapped
to the identity in the ball of radius $n$ of $F_d$ by the quotient map $\eta_n:F_d \to R_{n}$, among all
the proper limit (over free products) quotients of $R_{n-1}$ that admit an infinite descending chain of 
limit groups over free products that satisfy the conditions of  theorem 13.
If $R_{\infty}$ admits an infinite descending chain of limit groups over free products:
$$R_{\infty} \to L_1 \to L_2 \to \ldots$$
that satisfy the properties  in theorem 13, 
then the limit group (over free products) $L_1$ admits an infinite descending chain of limit
groups that satisfy the conditions of theorem 13,
and since it is a proper quotient of $R_{\infty}$, for large enough  index $n$, the quotient map
$\nu_n : F_d \to L_1$ maps to the identity strictly more elements of the ball of radius $n$ in the Cayley
graph of $F_d$, than the map $\eta_n:  F_d \to R_n$, a contradiction. Hence, $R_{\infty}$ does
not admit an infinite descending chain of limit groups over free products that satisfy the conditions of
theorem 13.

\smallskip
To continue the proof of theorem 13, i.e., to contradict the existence of the infinite descending
chain of limit groups over free products that satisfies the conditions  of the theorem,
we need a modification of the shortening
procedure that was used in [Se1] for ($F_k$) limit groups, and in [Se3] for limit groups over hyperbolic
groups.
Since
the description of the shortening procedure is rather long and involved, we prefer  not
to repeat it, and refer the interested reader to section 3 of [Se1]. 
The same construction that
appears in [Se1] applies to (strict) limit groups over free products.

Given a f.g.\ group $G$,
and a sequence of homomorphisms into free products: $\{u_s: G \to A_s*B_s\}$, that converges into a
(strict) limit group over free products, $L$, the shortening procedure constructs another (sub) 
sequence of homomorphisms
from a free group $F_d$ (where the f.g.\ group $G$ is generated by $d$ elements),
$\{v_{s_n}:F_d  \to A_{s_n}*B_{s_n}\}$, so that the sequence of homomorphisms $v_{s_n}$ converges to a 
limit group over free products
$SQ$, and there exists a natural epimorphism $L \to  SQ$, that maps the elliptic
elements in $L$, $E_{L}$, monomorphically into the elliptic elements in $SQ$,
$E_{SQ}$. 

\vglue 1pc
\proclaim{Definition 15} 
We call the 
limit group over free products, $SQ$, that is obtained by the shortening procedure,  
a shortening quotient of the limit group (over free products)
$L$. 
\endproclaim

By construction, a shortening quotient of a limit group over free products is, in particular, a quotient of that 
limit group. In the case of freely indecomposable $F_k$-limit groups, a shortening quotient is always a proper quotient
([Se1],5.3).
 If the limit group over free products that we start with, $L$, is strict, non-cyclic and 
admits no free decomposition in which the elements of $E_L$ can be conjugated into the factors,
a shortening quotient of $L$
is a proper quotient of it. More generally we have the following.

\vglue 1pc
\proclaim{Proposition 16} Let $G$ be a f.g.\ group, 
and let $\{u_s:G \to A_s*B_s\}$ be a sequence of  homomorphisms that converges into an action of a 
non-cyclic,  strict limit group over free products, $L$, on some real tree $Y$, where $L$ 
admits a (possibly trivial) free decomposition
in which the elliptic elements, $E_L$, can be conjugated into the factors, and so that there exists at least one
factor in this free decomposition, which is strict, non-cyclic, and freely indecomposable 
relative to its intersection with $E_L$.
Then every shortening quotient of 
$L$, obtained from the sequence $\{u_s\}$,  is a proper quotient of $L$ (in which non-trivial elliptic 
elements in $L$
are not mapped to non-trivial elliptic elements of the shortening quotient).
\endproclaim

\nfp 
Suppose that the f.g.\ group $G$ is generated by $d$ elements. A shortening quotient $SQ$ of $L$ is obtained
from a sequence of homomorphisms $\{v_{s_n}: F_d \to A_{s_n}*B_{s_n}\}$ that converges into $SQ$. Let $L_1$ be a
factor in a (possibly trivial) free decomposition of $L$, in which all the elements $E_L$ can be conjugated into the
various factors, so that the factor $L_1$ is a non-cyclic strict limit group (over free products), which is freely
indecomposable relative to its intersection with $E_L$.

Let $SQ_1$ be the image of $L_1$ in the shortening quotient $SQ$. Note that $SQ_1$ is a shortening quotient of $L_1$.  
By construction, the shortening quotient $SQ_1$ is a quotient of the non-cyclic, strict limit group over free
products $L_1$, which is freely indecomposable relative to its elliptic elements, $E_{L_1}$. If the sequence of
homomorphisms $\{v_{s_n}\}$, restricted to some f.g.\ preimage of $L_1$, 
has bounded stretching factors, i.e., if the shortening quotient $SQ_1$ is
not  strict, $SQ_1$ can not not be (entirely) elliptic, hence, it must be  freely decomposable or cyclic, so it
is a proper quotient of $L_1$. If $SQ_1$ is a strict limit group over free products,
 then the shortening
quotient $SQ_1$ is a proper
quotient of $L_1$ by the shortening argument that is used in the proof of claim 5.3 in [Se1]. 

\line{\hss$\qed$}

The shortening procedure, and proposition 16, enable us to obtain a $resolution$ of the limit group
$R_{\infty}$, with which we can associate a $completion$, into which $R_{\infty}$ embeds. This completion
enables us to present $R_{\infty}$ as a f.g.\ group which is finitely presented over some of its elliptic
subgroups. Since  theorem 13 assumes that the successive maps along the infinite descending chains
under consideration, $\tau_i$, do not map non-trivial elliptic elements to the identity
element, it is implied that elliptic subgroups embed along the sequences
under consideration. This implies that for large enough $n$, $R_n$ is a proper quotient of $R_{\infty}$,
which contradicts the fact that $R_{\infty}$ is a proper quotient of all
the limit groups (over free products), $\{R_n\}$, that appear in the infinite descending chain we constructed.

\vglue 1pc
\proclaim{Proposition 17} Let $R_{\infty}$ be the direct limit of the sequence of limit groups
over free products that we constructed (in order to prove theorem 13), $\{R_n\}$.
Then there exists a finite sequence of limit groups over free products:
$$R_{\infty}  \to L_1 \to L_2 \to \ldots \to L_s$$
for which:
\roster
\item"{(i)}" $L_1$ is a shortening quotient of $R_{\infty}$, and $L_{i+1}$ is a shortening
quotient of $L_i$, for $i=1,\ldots,s-1$.

\item"{(ii)}" The epimorphisms along the sequence are proper epimorphisms, and non-trivial elliptic elements 
in $L_i$  
are mapped to non-trivial elliptic elements in  $L_{i+1}$.

\item"{(iii)}" $L_s=H_1* \ldots *H_r*F_t$ where the factors, $H_1,\ldots,H_r$, are elliptic, and the entire
elliptic set, $E_{L_s}$, is the union of the conjugates of $H_1,\ldots,H_r$. $F_t$ is a (possibly trivial) free group. 

\item"{(iv)}" The resolution:
$R_{\infty}  \to L_1 \to L_2 \to \ldots \to L_s$ is a strict resolution ([Se1],5), i.e., in each level
non-QH, non-virtually-abelian
vertex groups in the virtually abelian JSJ decomposition are mapped monomorphically into the limit group in the next
level, and QH vertex groups are mapped into non-virtually-abelian, non-elliptic subgroups. 

\item"{(v)}" The 
constructed resolution is well-structured (see definition 1.11 in [Se2] for a well-structured resolution). 

\endroster
\endproclaim

\nfp By lemma 14 and proposition 16 a shortening quotient of $R_{\infty}$ is a proper quotient of it. 
Furthermore, non-trivial elliptic elements
in $R_{\infty}$ are mapped to non-trivial elliptic elements in the  shortening quotient.   
Hence, we set $L_1$ to be a shortening quotient of $R_{\infty}$. If from the sequence of (shortened) homomorphisms 
that was used to construct $L_1$, it's possible to extract a subsequence that satisfy the properties
of lemma 14, we continue with this subsequence, and use it to get a shortening
quotient $L_2$ of $L_1$, which by proposition 16, is a proper quotient of $L_1$. Continuing this
process iteratively, and recalling that every descending chain of limit groups over free products
that starts with $R_{\infty}$ and satisfies the assumptions of the statement of theorem 13, terminates after finitely
many steps,  
 we finally get the sequence of proper epimorphisms:
$$R_{\infty}  \to L_1 \to L_2 \to \ldots \to L_s.$$

Parts (i) and (ii) follow immediately from the construction of the descending finite sequence of shortening
quotients, and part (iii) follows, since by lemma 14 and proposition 16, the descending sequence of
shortening quotients terminates, precisely when the obtained
limit group is the free product of elliptic factors and a free group. Part (iv) follows since each shortening
quotient in the sequence is constructed from homomorphisms that converge into the previous limit group
in the sequence that were modified by modular automorphisms. Part(v) follows since like in the case of free
and hyperbolic groups, every strict Makanin-Razborov resolution, i.e., a resolution that is obtained from
a sequence of shortening quotients, is well-structured (see definition 1.11 in [Se2]).

\line{\hss$\qed$}

Proposition 17 constructs from a subsequence of the homomorphisms, $\{h_n:F_d \to A_n*B_n\}$, a
well-structured resolution of the limit group over free products, $R_{\infty}$, that terminates in a limit group
$L_s$ which is a free product of elliptic subgroups and a (possibly trivial) free group. In section 1 of [Se2],
a $completion$ is constructed from  a given  well-structured resolution (see definition 1.11 in [Se2]
for a well-structured resolution). This
construction that generalizes in a straightforward way to well-structured resolutions of limit groups over
torsion-free hyperbolic groups in [Se3], generalizes in a straightforward way to well-structured
resolutions of limit groups over free products. For the detailed construction of the completion
see definition 1.12 in [Se2].

We denote the completion of the well-structured resolution that is constructed in proposition 17,
$Comp(Res)$.
By definition 1.12 and lemma 1.13 in [Se2], each of the limit groups, $R_{\infty},L_1,\ldots,L_s$ is
embedded into the completion of the constructed resolution, $Comp(Res)$. All the (virtually
abelian) edge groups that are not connected to boundary elements of QH vertex groups (edge groups
that are connected to QH vertex groups  are always  cyclic), 
and all the virtually abelian vertex groups, contain abelian subgroups as subgroups of index at most 2.
Furthermore, these abelian groups are non-elliptic subgroups of the associated limit groups, 
$R_{\infty},L_1,\ldots,L_s$. Since the only non-elliptic abelian subgroups of the terminal limit
group, $L_s$, are infinite cyclic, all the edge groups, and all the vertex groups, 
that appear in all the levels of the completion, $Comp(Res)$, are finitely generated. In particular all
the edge groups and all the vertex groups
that appear in the virtually abelian JSJ decompositions of the limit groups over free
products, $R_{\infty},L_1,\ldots,L_s$, are finitely generated.

Let $\rho: R_{\infty} \to Comp(Res)$ be the embedding of the limit group over free products, $R_{\infty}$,
into the completion of the constructed resolution: $R_{\infty} \to L_1 \to \ldots \to L_s$. $\rho(R_{\infty})$
being a f.g.\ subgroup of $Comp(Res)$ inherits a (finite) 
virtually abelian decomposition from the virtually abelian decomposition
that is associated with the top level of $Comp(Res)$. The edge groups in that 
inherited (finite) virtually abelian decomposition
are  subgroups of f.g.\ virtually abelian groups, hence, f.g.\ virtually abelian groups.
The vertex groups in that virtually abelian decomposition can either be conjugated into subgroups of 
a lower level of
the completion,  or they can be conjugated into QH groups or into f.g.\ virtually abelian groups. f.g.\ subgroups
of virtually abelian groups are again f.g.\ virtually abelian. f.g.\ subgroups of Fuchsian groups are free products of
f.g.\
Fuchsian 
groups and  f.g.\ virtually free groups. Hence $R_{\infty}$ can be reconstructed from finitely many f.g.\ groups
that can be conjugated into lower level of the completion, $Comp(Res)$, and finitely many f.g.\ Fuchsian groups,
 f.g.\ virtually free groups, and f.g.\ virtually abelian groups, by performing free products and 
free products with amalgamation and HNN extensions along f.g.\ virtually abelian groups.

Continuing with this decomposition procedure along the lower levels of the completion, $Comp(Res)$, we get
that the subgroup $\rho(R_{\infty})$ (that is isomorphic to $R_{\infty}$) can be reconstructed from finitely many
f.g.\ elliptic subgroups in $Comp(R_{\infty})$, 
and finitely many f.g.\ Fuchsian groups,
 f.g.\ virtually free groups, and f.g.\ virtually abelian groups, by performing finitely
many operations of free products and 
free products with amalgamation and HNN extensions along f.g.\ virtually abelian groups. 
In particular, $R_{\infty}$ is obtained
from finitely many elliptic subgroups of $R_{\infty}$ by adding finitely many generators and relations.

By construction, the limit group (over free products), $R_{\infty}$, 
is the direct limit of the decreasing sequence of
limit groups, $\{R_n\}$, which are all quotients of some free group, $F_d$. Every f.g.\ subgroup of
$F_d$ that is mapped to an elliptic subgroup in $R_{\infty}$, is mapped to elliptic subgroups in
all the  limit groups, $R_n$, for all $n>n_0$ for some index $n_0$. 
$R_{\infty}$ is generated by
finitely many f.g.\ elliptic subgroups and finitely many virtually abelian, virtually free, and f.g.\
Fuchsian groups together with finitely many Bass-Serre generators that are added in each of the performed
HNN extensions (along f.g.\ virtually abelian subgroups). 
Since these last groups are all f.p.\ and elliptic subgroups in each of the limit groups $R_n$
are mapped monomorphically into $R_{\infty}$ by our assumptions on the decreasing sequence, $\{R_n\}$, 
there exists some index $n_1$, so that for all $n>n_1$, the limit group $R_n$ is generated by 
finitely many elements that are mapped to the Bass-Serre elements that are used in constructing $R_{\infty}$,
finitely many subgroups that are isomorphic to the f.g.\ virtually
abelian, f.g.\ virtually free, and f.g.\ Fuchsian groups, and finitely many
elliptic
subgroups that are isomorphic to the  f.g.\ 
elliptic subgroups that altogether generate $R_{\infty}$. Since $R_{\infty}$  is generated by these subgroups
and the Bass-Serre elements by imposing finitely many relations, there exists some index $n_2$, so that for every $n>n_2$
these relations hold in $R_n$, which implies that $R_n$ is a quotient of $R_{\infty}$ using a quotient map
that maps the fixed generating set of $R_{\infty}$ (the images of a fixed basis of $F_d$) to the fixed
generating set of $R_n$. This implies that $R_{n+1}$ is a proper quotient of $R_{\infty}$ by a quotient map
that maps the fixed generating set of $R_{\infty}$ to the fixed generating set of $R_{n+1}$, which 
contradicts the construction of $R_{\infty}$ as the direct limit of the decreasing sequence of limit groups
over free products, $\{R_n\}$. This finally implies the d.c.c.\ that is stated in theorem 13.

\line{\hss$\qed$}

Theorem 13 proves a basic d.c.c.\ that holds for limit groups over free products. This d.c.c.\ is
weaker than the ones proved for limit groups over free and hyperbolic groups ([Se1],[Se3]).
Indeed it is stated only for decreasing sequences of limit groups over free products for which the
successive maps do not map non-trivial elliptic elements to the identity. Still, this d.c.c.\
is the basis for our analysis of limit groups over free products, and for the analysis of solutions to systems
of equations over free products.

We start with the following theorem, which is a rather immediate corollary of the d.c.c.\ that is stated in theorem
13, 
that associates a resolution with a given limit group over free products, a resolution that has similar properties
to the resolution described in proposition 17.

\vglue 1pc
\proclaim{Theorem 18} Let $L$ be a limit group over free products.
Then there exists a finite sequence of limit groups over free products:
$$L  \to L_1 \to L_2 \to \ldots \to L_s$$
for which:
\roster
\item"{(i)}" $L_1$ is a shortening quotient of $L$, and $L_{i+1}$ is a shortening
quotient of $L_i$, for $i=1,\ldots,s-1$. In particular, elliptic elements in $L_i$ are mapped monomorphically 
to elliptic elements in $L_{i+1}$.

\item"{(ii)}" The epimorphisms along the sequence are proper epimorphisms.

\item"{(iii)}" $L_s=H_1* \ldots *H_r*F_t$ where the factors, $H_1,\ldots,H_r$, are elliptic, and the entire
elliptic set, $E_{L_s}$, is the union of the conjugates of $H_1,\ldots,H_r$. $F_t$ is a (possibly trivial) free group. 

\item"{(iv)}" The resolution:
$L  \to L_1 \to L_2 \to \ldots \to L_s$ is a strict resolution ([Se1],5), i.e., in each level
non-QH, non-virtually-abelian
vertex groups and edge groups 
in the virtually abelian JSJ decomposition are mapped monomorphically into the limit group in the next
level, and QH vertex groups are mapped into non-virtually-abelian, non-elliptic subgroups. 

\item"{(v)}" The 
constructed resolution is well-structured (see definition 1.11 in [Se2] for a well-structured resolution). 
As a corollary, the limit group (over free products) $L$ is embedded into the completion of the well-structured
resolution:
$$L  \to L_1 \to L_2 \to \ldots \to L_s$$ 
so that all the elliptic elements in $L$ are mapped into conjugates of the elliptic subgroups, $H_1,\ldots,H_r$,
of $L_s$.
\endroster
\endproclaim

\nfp Theorem 18 generalizes the resolution that was constructed for the limit group (over free products)
$R_{\infty}$, to general limit groups over free products. To prove proposition 17 we used the d.c.c.\ for
resolutions of $R_{\infty}$ for which the epimorphisms that are associated with them 
do not map non-trivial elliptic elements to the identity element, 
that follows from the construction of $R_{\infty}$. Theorem 13 proves that such a d.c.c.\ holds for 
resolutions of an arbitrary limit group over free products, for which the associated epimorphisms do
not map non-trivial elliptic elements to the identity element.
With this general d.c.c.\ the proof of proposition 17 generalizes to
general limit groups over free products.

\line{\hss$\qed$}

\vglue 1.5pc
\centerline{\bf{\S3. Finitely Presented Groups}}
\medskip

Theorem 13 proves the basic d.c.c.\ for limit groups over free products, and theorem 18 associates
a resolution with each such limit group, hence, it embeds each limit group over free products into a
completion, where this completion is a tower over a limit group which is a free product of a (possibly trivial)
free group with a (possibly empty) finite collection of f.g.\ elliptic subgroups. 

When considering limit groups over free products we analyzed sequences of homomorphisms from a f.g.\
group into free products. Since our goal is to obtain a structure theory for sets of solutions to
systems of equations, and the group that is associated formally with a finite system of equations is
f.p.\ and not only f.g.\ we may assume that the limit groups over free products that we are considering
are obtained from sequences of homomorphisms from a f.p.\ group into free products (and not only
from a f.g.\ one).

As we will see in the sequel, if we attempt to construct a Makanin-Razborov diagram that is associated
with a f.p.\ group, we will need to consider only f.g.\ limit groups over free products that are recursively 
presented, i.e., limit groups that can be embedded into f.p.\ groups. A modification or a strengthening of
the existence of such an embedding is a 
key for obtaining further d.c.c.\ that will eventually allow the construction of a Makanin-Razborov
diagram over free products for a given f.p.\ group. We start with the
following simple observation.

\vglue 1pc
\proclaim{Proposition 19} Let $G$ be a f.p.\ group, and let $L$ be a limit group over free products which is
a  quotient of $G$. Then there  exists a limit group over free products $\hat L$ with the following properties:
\roster
\item"{(1)}" there is a f.p.\ completion, $Comp$, which is a tower  over a free product of finitely many 
f.p.\ elliptic subgroups and a free group, so that $\hat L$ embeds into $Comp$, and the elliptic elements in $\hat L$
are mapped into conjugates of the finitely many elliptic factors in the free decomposition that is associated with
the limit group that appears in the terminal level of the completion $Comp$. 

\item"{(2)}" either $\hat L=L$ or $\hat L > L$ (see definition 12 for the relation $>$ on limit groups over 
free products).
\endroster
\endproclaim

\nfp By theorem 18, the limit group (over free products) $L$ admits a well-structured resolution:
$$L  \to L_1 \to L_2 \to \ldots \to L_s$$ 
and $L_s$ admits a free product decomposition:
$L_s=H_1* \ldots *H_r*F_t$ where the factors, $H_1,\ldots,H_r$, are elliptic, and the entire
elliptic set, $E_{L_s}$, is the union of the conjugates of $H_1,\ldots,H_r$. $F_t$ is a (possibly trivial) free group. 

Furthermore, with this resolution it is possible to associate a completion, $Comp_1$, and the limit group $L$ embeds into 
this completion, 
so that all the elliptic elements in $L$ are mapped into conjugates of the elliptic subgroups, $H_1,\ldots,H_r$,
of $L_s$ (the groups $L_1,\ldots,L_s$ admit natural embeddings into the various levels of the completion, 
$Comp_1$, and
the elliptics in each of these limit groups are mapped into conjugates of $H_1,\ldots,H_r$ in the completion
$Comp_1$).

Since $L$ is embedded into the completion $Comp_1$, $G$ is naturally mapped into $Comp_1$. By construction,
the completion $Comp_1$ is
built as a tower over the terminal limit group $L_s$. If $Comp_1$ is f.p.\ we obtained the conclusion of the
proposition, as we can take $\hat L=L$, and $\hat L$ is embedded into the f.p.\ completion $Comp_1$. 
Hence, we may assume that $Comp_1$ is not finitely presented, i.e., at least one of 
the factors, $H_1,\ldots,H_r$, is not finitely presented.
In that case we gradually replace $Comp_1$ by a f.p.\ completion into which
$G$ is mapped.

Each of the factors of $L_s$, $H_1,\ldots,H_r$,  is f.g.\ so it is a quotient of some f.g.\ free group. 
Let $F^1,\ldots,F^r$ be f.g.\ free groups that $H_1,\ldots,H_r$ are quotients of. We
start the construction of a f.p.\ completion that replaces the completion $Comp_1$,  
with a tower $T_2$ that has in 
its base level the free group $F^1*\ldots*F^r*F_t$, and the next (upper) levels are connected to the lower levels
of the constructed tower,
precisely as they are connected in the completion, $Comp_1$, i.e., using the same graphs of groups, just that the 
 group that is associated with the lowest level in $Comp_1$, which is $L_s=H_1*\ldots*H_r*F_t$, 
is replaced by the free group, $F^1*\ldots*F^r*F_t$.

Note that $T_2$ is a tower, but it is not necessarily a completion (see definition 1.12 in [Se2]), 
as in general there are no retractions
from a group that is associated with a certain level onto the group that is associated with the level below it.
Each of the levels above the base level in $T_2$ is  constructed using a (finite) graph of groups, in which
some vertex groups are the groups that are associated with the lower level in $T_2$. Hence, the group that is
associated with a level above the base level,  
is obtained from a free product of the group that is associated with the lower level with a f.p.\ group
by imposing finitely many relations. Furthermore, the graphs of groups that are associated with the different
levels in $T_2$ are similar to the  graphs of groups that are associated with the corresponding levels in the
completion $Comp_1$, and differ from $Comp_1$ only in the groups that are associated with the base level. 

Since each of the groups that are associated with the upper levels in $T_2$ is obtained from a free 
product with a f.p.\ group by imposing finitely many relations, and since the graphs of groups that are
associated with the upper levels have similar structure as the corresponding graphs of groups that
are associated with the levels of the completion $Comp_1$, and these graphs of groups differ only in the structure
of the group that is associated with the base level, it is enough to impose only finitely many relations 
from the defining relations of the various factors of the limit groups that is associated with
the base level of $Comp_1$, $L_s$,
$H_1,\ldots,H_r$, on  the associated free groups, $F^1,\ldots,F^r$, so that if we replace the group that
is associated with the base level of $T_2$,
$F^1*\ldots,F^r*F_t$, with the obtained f.p.\ quotient, $V_1*\ldots*V_r*F_t$, and construct from it a tower,
$T_3$, by imitating the construction of $Comp_1$ and $T_2$ (i.e., with similar graphs of groups in all the upper
levels), $T_3$ will be a completion.
     
$T_3$ is a completion, but it may be that the f.p.\ group $G$ is not mapped into it. 
$G$ is mapped into the completion $Comp_1$. Hence, once again, since $G$ is f.p.\ it is enough to impose 
only finitely many relations from the defining relations of the various factors, $H_1,\ldots,H_r$, on
the factors, $V_1,\ldots,V_r$, so that if we replace the group that is associated with the base level in $T_3$
with the obtained f.p.\ quotient, $M_1*\ldots*M_r*F_t$, and construct from it a tower $T_4$ by imitating the 
construction of the towers $Comp_1$, $T_2$, and $T_3$, $T_4$ is a f.p.\ completion, and $G$ maps into it. 

Furthermore, the map from $G$ into the completion $Comp_1$, is a composition of the maps from $G$ to $T_4$, 
composed with the natural quotient map from $T_4$ to $Comp_1$. Hence, if we denote the image of $G$ in $T_4$
by $\hat L$, then $\hat L$ is a limit group over free products, its set of elliptics is precisely the 
intersection of $\hat L$ with the set of conjugates of $M_1,\ldots,M_r$, and either $\hat L$ is
isomorphic to $L$ and the natural isomorphism from $\hat L$ onto $L$ maps the elliptics in $\hat L$ 
monomorphically
onto the elliptics in $L$, or the natural epimorphism from $\hat L$ onto $L$ has a non-trivial kernel,
and this epimorphism maps the elliptics in $\hat L$ onto the elliptics in $L$, in which case $\hat L > L$. 

\line{\hss$\qed$}

Proposition 19, the d.c.c.\ proved in theorem 13, and the resolution that is associated with a limit group
over free products in theorem 18, enable us to prove 
that there are maximal elements in the set of all limit groups over free products that
are all quotients of a (fixed) f.p.\ group $G$,
and that there are only finitely many equivalence classes of such maximal elements. The existence of maximal 
elements in the set of limit quotients is valid even for f.g.\ groups.

\vglue 1pc
\proclaim{Proposition 20} Let $G$ be a f.g.\ group. 
Let $R_1,R_2,\ldots$ be a sequence of limit groups over free products that are all quotients
of the f.g.\ group $G$, and for which:
 $$R_1 \, < \, R_2 \, < \, \ldots$$ 
 Then there exists a limit group over free products $R$ that is a quotient of $G$,
so that for every index $m$, $R > R_m$.
\endproclaim

\nfp Identical to the proof in the free and hyperbolic groups cases (see proposition 1.20 in [Se3]). 

\line{\hss$\qed$}

Proposition 20 proves that given an ascending chain of limit quotients (over free products) of a f.g.\ 
group $G$, there exists a limit quotient of $G$ that bounds all the limit groups in the sequence.
Hence, we can apply Zorn's lemma (it is enough to consider countable ascending chains in case of quotients
of a f.g.\ group), and obtain 
maximal limit quotients (over free products) of any given f.g.\ group, and every limit
quotient of a f.g.\ group is dominated by a maximal limit quotient of that group.

Proposition 19 proves that if $G$ is in addition f.p.\ then if $R$ is a limit 
quotient of $G$ (over free products), then there exists a limit group over free products $L$, 
that is either isomorphic to $R$
or $L>R$, so that $L$ embeds in a f.p.\ completion. Hence, if we are interested in maximal limit quotients
(over free products) of a f.p.\ group $G$, it is enough to consider limit quotients of $G$ that
embed in f.p.\ completions, and there are clearly at most countably many such limit quotients.


In case a group $G$ is f.p.\ the existence of maximal limit quotients, 
and the existence of an embedding of maximal limit quotients of a f.p.\ group $G$
into f.p.\ completions, imply the finiteness of
the (equivalence classes of) maximal limit quotients (over free products) of a f.p.\ group.  

\vglue 1pc
\proclaim{Theorem 21} Let $G$ be a f.p.\ group. 
Then there are only finitely many equivalence classes of maximal elements in 
the set of limit quotients (over free products)  of $G$, and each of these maximal elements embeds in a f.p.\
completion.
\endproclaim

\nfp  Let $G$ be a f.p.\ group. Since all its maximal limit quotients over free products can be embedded
into f.p.\ completions, there are at most  countably many maximal limit quotients of $G$ (over free products).
Suppose that there are infinitely many non-equivalent maximal limit quotients of $G$, and
let $R_1,R_2,\ldots$  be the  infinite sequence of
(non-equivalent) maximal limit quotients (over free products) of $G$. Each $R_i$ is
equipped with a given quotient map $\eta_i:G \to R_i$, hence, fixing a
generating set for $G$, we fix a generating set in each of the $R_i$'s. 
i.e.,
we have maps $\nu_i:F_d \to R_i$ (assuming $G$ has rank $d$), that factor through the epimorphism $F_d \to G$.

For each index $i$ we look at the collection of words of length 1 in $F_d$ that
are mapped to the identity,  and those that are mapped to elliptic elements by $\nu_i$. There is a subsequence of
the $R_i$'s for which this (finite) collection of words is identical.
Starting with this subsequence, for each $R_i$ (from the subsequence) we look at the collection of
words of length 2 in $F_d$ that are mapped to the identity and those that are mapped
to elliptic elements by $\nu_i$, and
again there is a subsequence for which this (finite) collection is
identical. We continue with this process for all lengths $\ell$ of words
in $F_d$, and look at the diagonal sequence (that we  denote
$R_{i_1},R_{i_2},\ldots$).

We choose homomorphisms $h_j: F_d \to A_j*B_j$, that factor through the map $F_d \to G$,
 so that for words $w$
of length at most $j$, $h_j(w)=1$ iff $\nu_{i_j}(w)=1$, and $h_j(w)$ is elliptic iff $\nu_{i_j}(w)$
is elliptic
(we can choose such homomorphisms  since
$R_{i_j}$ is a limit quotient of $G$). After passing to a subsequence, the homomorphisms
$h_{j}$ converge into a limit group over free products $M$, which is a limit quotient of $G$. Note that in the
(canonical) map $F_d \to M$, the elements of length at most $j$ that are
mapped to the identity, and those that are mapped to be elliptic,
 are precisely those that are mapped to the identity and those that are mapped to be elliptic
by the map $\nu_{i_j}: F_d \to R_{i_j}$.

$R_1,R_2,\ldots$ form  the entire list of maximal limit quotients of $G$ over free products. 
We construct a new sequence 
of homomorphisms: $f_j : F_d \to C_j*D_j$ that factor through the quotient map
$F_d \to G$. First, $f_j$ has the same property as $h_j$, i.e., the elements
of length at most $j$ that are mapped to the identity by $f_j$ are
precisely those that are mapped to the identity by $\nu_{i_j}: F_d \to R_{i_j}$, and the
elements of length at most $j$ that are mapped to be elliptic by $f_{j}$ are
precisely those that are mapped to be elliptic by $\nu_{i_j}: F_d \to R_{i_j}$. Second,
since  $R_{i_j}$ is maximal and is not equivalent to $R_1,\ldots,R_{i_j-1}$, there must exist some elements
$u_1,\ldots,u_{i_j-1} \in F_d$ so that for each index $s$, $1 \leq s \leq i_j-1$, either $u_s$ 
is mapped to the identity in $R_s$, but $u_s$ is mapped
to a non-trivial element in $R_{i_j}$ by $\nu_{i_j}$, or $u_s$ is mapped to an elliptic element in $R_s$,
but $u_s$ is mapped to a non-elliptic element in $R_{i_j}$ by $\nu_{i_j}$. If the first holds,
we require that $f_j(u_s) \neq 1$, and if the second holds we require that $f_j(u_s)$ is not elliptic.

The sequence of homomorphisms, $\{f_j\}$, converges into the limit group (over free products) $M$.
We look at a subsequence of the homo. $\{f_j\}$, so that the subsequence and its shortenings converge into 
a resolution of $M$ that satisfy
the properties that are listed in theorem 18, $M \to L_1 \to L_2 \to \ldots \to L_s$
(we still denote the subsequence  $\{f_j\}$). 

With the resolution $M \to L_1 \to L_2 \to \ldots \to L_s$, which is a well-structured resolution by
construction, we can naturally associate a completion. Let $Comp$ be this completion. Since $G$ is naturally
mapped onto the limit group $M$, there exists a natural map, $\rho: G \to Comp$, that factors through the map
$G \to M$. 

\noindent
Note that by construction, the completion $Comp$ is obtained from the terminal limit group, $L_s$,
of the given resolution of $M$, by adding finitely many generators and relations. Since the group $G$ is
f.p.\ we can repeat the argument that was used to prove proposition 19, and replace the terminal limit group
$L_s$ with a (possibly the same) f.p.\ group $L^1_s$ that maps onto $L_s$, and starting with $L^1_s$ construct
a completion, $Comp^1$, that has the same structure as the completion $Comp$, except that the terminal
limit group (over free products) of $Comp^1$ is $L^1_s$, whereas the terminal limit group of the completion
$Comp$ is $L_s$. Furthermore, the group $G$ maps into $Comp^1$, and since $Comp^1$ is finitely presented, 
there exists a subsequence of the 
sequence of homomorphisms $\{f_j\}$, that factor through the completion $Comp^1$.

Let $M^1$ be the image of $G$ in $Comp^1$. $M^1$ is a limit quotient of $G$  (over free products), so there must exist some
maximal limit quotient of $G$, that we denote $R_b$, so that $R_b$ is either equivalent to $M^1$ or $R_b>M^1$.
Now, there exists a subsequence of the homomorphisms $\{f_j\}$ that factor through the limit group
$M^1$, hence, this subsequence of homomorphisms factor through the maximal limit quotient $R_b$. By construction,
each of the homomorphisms $f_j$ does not factor through any of the maximal limit groups, $R_1,\ldots,R_{i_j-1}$. 
Hence, for large enough $j$, none of the homomorphisms $f_j$ factor through the maximal limit quotient $R_b$,
a contradiction. Therefore, $G$ admits only finitely many maximal limit quotients (over free products),
and by proposition 19, each of the maximal limit quotients of $G$ embeds into a f.p.\ completion.

\line{\hss$\qed$}

Theorem 21 proves the existence of finitely many limit quotients of a given f.p.\ group. Hence,
it gives the first level of a Makanin-Razborov diagram of a f.p.\ group over free products, and it proves that
the groups that appear in the first level of the Makanin-Razborov diagram of a f.p.\ group over free 
products are canonical (i.e., they are an invariant of the f.p.\ group). Still, the construction of maximal
limit groups over free products, and the proof that there are only finitely many (equivalence classes of) maximal
quotients of a f.p.\ group (over free products), does not generalize in a straightforward way to allow the
 construction of the next levels in the Makanin-Razborov diagram. 

Furthermore, 
theorem 13
proves the basic d.c.c.\ that is required for analyzing limit groups over free products. However, it is
not sufficient for the construction of a Makanin-Razborov diagram as it guarantees the termination
of strict resolutions, but not of general resolutions in the diagram (if we try to imitate
the construction over free and hyperbolic groups). Hence, to construct a finite 
Makanin-Razborov diagram we will need to construct the next levels in the diagram, and in addition to
prove an additional d.c.c.\ that will guarantee the termination of the construction after finitely many steps.

Let $G$ be a f.p.\ group. 
We start the construction of the Makanin-Razborov diagram over free products of $G$ with the finite
collection of (equivalence classes of) maximal limit quotients of $G$, according to theorem 21. We
continue by studying the  homomorphisms of each of the maximal limit
quotients of $G$ into free products. 
As in the construction of Makanin-Razborov diagrams over free and hyperbolic groups, we 
continue by modifying (shortening) these  homomorphisms using the modular groups that are associated with
the maximal limit quotients (over free products) of the given f.p.\ group $G$.

Let $L$ be one of the maximal limit quotients (over free products) of $G$, and let $E_L$ be its set of elliptics.
First, we factor $L$ into its most refined free decomposition in which the elements in $E_L$ are
elliptic (i.e., contained in conjugates of the factors), $L=U_1*\ldots*U_m*F_t$, where $F_t$ is a (possibly
trivial free group, and the elements in $E_L$ can be conjugated into the various factors, $U_1,\ldots,U_m$.

$(L,E_L)$ is a (maximal) limit quotient of $G$ (over free products), hence, $L$ is obtained as a limit
of a sequence of homomorphisms $\{h_n:G \to A_n*B_n\}$. $G$ is f.p.\ and is  mapped onto $L$, and
$L$  admits the free decomposition, $L=U_1*\ldots*U_m*F_t$, where the elliptic elements in $E_L$ can
be embedded into the various factors $U_1,\ldots,U_m$. By the argument that is used to prove proposition 19, 
there exist
finitely presented groups $M_1,\ldots,M_m$ so that the map $G \to L$ factors as:
$$G  \ \to \ M_1*\ldots*M_m*F_t \ \to \ U_1*\ldots*U_m*F_t$$
where the two maps are onto, and for each index $i$, $1 \leq i \leq m$, $M_i$ is mapped onto $U_i$. Since  
the sequence of homomorphisms $\{h_n\}$ of $G$ 
converges into $(L,E_L)$, and the group $M_1*\ldots*M_m*F_t$ is f.p.\ and the map from $G$ to $L$ factors through it,
for large enough $n$ the homomorphisms $\{h_n\}$ factor through the map
$G \to M_1*\ldots*M_m*F_t$.  Now, if we apply the proof of proposition 19, it follows that there are
$m$ f.p.\ completions (over free products), $Comp_1,\ldots,Comp_m$, so that each of the factors $U_i$ is embedded
into the completion $Comp_i$ so that the elliptics in $U_i$ are mapped into elliptics in $Comp_i$
(and only elliptics in $U_i$ are mapped into elliptics in $Comp_i$), and
there exist  maps:
$$M_1*\ldots*M_m*F_t \ \to \ U_1*\ldots* U_m*F_t \ \to \ Comp_1*\ldots*Comp_m*F_t $$
that extend the embeddings from $U_i$ to $Comp_i$, for $1 \leq i \leq m$.

Hence, we may continue with each of the factors $U_i$ of $L$ in parallel. $U_i$ is a maximal limit quotient
(over free products) of the f.p.\ group $M_i$, and by proposition 21 it is embedded into a f.p.\ completion $Comp_i$.

\noindent
Therefore, we may assume that in the sequel,  we are given a f.p.\ group $G$, and a maximal limit quotient of
it, that we still denote, $(L,E_L)$, and the limit quotient $L$ is freely indecomposable relative to the elliptic
subset $E_L$ (i.e., $L$ admits no non-trivial free decomposition in which the elements in $E_L$ can
be conjugated into the factors).

\noindent
With (the factor) $L$ and $E_L$ we naturally associate its virtually abelian JSJ decomposition over
free products (theorem 11). We also associate with $(L,E_L)$ the collection of homomorphisms of $G$ 
into free products that factor through $(L,E_L)$.

Fixing a (finite) generating set of a limit group (over free products) $L$, and given a 
homomorphism $h: L \to A*B$, we look at a shortest homomorphism
among those that are obtained by precomposing $h$ with a modular automorphism of $L$ that is contained
in the modular group of automorphisms of $L$ that is associated with the virtually abelian JSJ decomposition over
free products of $L$
(relative to $E_L$). A limit group over free products that is the limit of a sequence of such shortest 
homomorphisms is called a shortening quotient, and denoted $SQ$. Note that this definition of a shortening
quotient is different than the more restricted 
one given in definition 15, as in particular, the natural map from a limit
group over free products, $L$, onto a shortening quotient $SQ$ of $L$, is not always monomorphic on the set
of elliptic elements in $L$, $E_L$. Still, 
like in proposition  16, if a shortening 
quotient is 
not elliptic it is a proper quotient of the limit group $L$.

\vglue 1pc
\proclaim{Lemma 22 (cf. proposition 16)} Let $L$ be a limit group over free 
products, and let $E_L$ be its set of elliptics. Suppose that $L$ admits no free decompositions in
which the elements in $E_L$ can be conjugated into the factors. Then every shortening quotient of $L$ which is
not (entirely) elliptic
is a proper quotient of it.
\endproclaim

\nfp Identical to the proof of proposition 16.

\line{\hss$\qed$}

Like limit quotients (over free products) of a f.g.\ group, every ascending sequence of shortening 
quotient of a limit group over free products is bounded by a shortening quotient of that limit group.

\vglue 1pc
\proclaim{Lemma 23} Let $L$ be a f.g.\ limit group over free products. 
Let $SQ_1,SQ_2,\ldots$ be a sequence of shortening quotients of $L$,
for which:
 $$SQ_1 \, < \, SQ_2 \, < \, \ldots$$ 
 Then there exists a shortening quotient $SQ$ of $L$,
so that for every index $m$, $SQ > SQ_m$.
\endproclaim

\nfp Identical to the proof in the hyperbolic group case (proposition 1.20 in [Se3]).

\line{\hss$\qed$}

By Zorn's lemma and lemma 23 it follows that there are maximal elements in the set of shortening quotients
of a f.g.\ limit group over free products. We call such a maximal element, a $maximal$ $shortening$
$quotient$. By lemma 22, if the limit group (over free products) $L$ does
not admit a free product in which the elliptic elements in $L$, $E_L$, can be conjugated into the factors,
every maximal shortening quotient of $L$ that is not entirely elliptic is a proper quotient of $L$.

\vglue 1.5pc
\centerline{\bf{\S4. Covers of Limit Quotients and their Resolutions}}
\medskip

The first level in the Makanin-Razborov diagram over free products 
of a f.p.\ group $G$ consists of the finitely many maximal 
limit quotients of $G$ (theorem 21). Over free and hyperbolic groups we continued to the next level in
the diagram by proving that there are only finitely many (equivalence classes of) maximal shortening quotients.
Over free product we need to prove a finiteness result for shortening quotients and their (strict)
resolutions, that will enable
us to continue to the next level, and so that the next levels will
 be constructed in a way for which a termination can be proved.

In order to prove that there are only finitely many maximal limit quotients over free products
of a f.p.\ group over
free products (theorem 21), we first showed that any maximal limit quotient can be embedded into a f.p.\
completion (proposition 19). 
For maximal shortening quotients of a f.g.\ limit group over free products we were not able
to prove such a statement. For the continuation of the diagram, we first prove an observation
that holds for all  the (proper) limit quotients of a given limit group over free products, 
that still allows us to construct the 
Makanin-Razborov diagram over free products for a f.p.\ group, although we loose some of the canonical
properties of the diagrams over free and hyperbolic groups.

Given a limit group over free products, $L$, and a limit quotient $M$ of $L$, theorem 24 associates a cover, $CM$,
with $M$. $CM$ is a limit quotient of $L$, if $L> M$, then $L>CM$ and $M$ is a limit quotient of $CM$.
The main property of the cover $CM$ that is used in the sequel (and is not always true for $M$) is that
$CM$ can be embedded into a completion, $Comp_{CM}$, and $Comp_{CM}$ is f.p.\ relative to the elliptic
subgroups of the given limit group $L$, i.e., $Comp_{CM}$ is generated from the elliptic subgroups in $L$ by
adding finitely many generators and relations. In particular, this implies that if $L$ is recursively
presented so is the cover $CM$.

\vglue 1pc
\proclaim{Theorem 24} Let $L$ be a f.g.\ limit group over free products, 
and let $E_L$ be its
set of elliptics. Let $M$ be a limit quotient of $L$ (over free products), with set of elliptics, $E_M$,
and with a quotient map,
$\eta: L \to M$ that maps $E_L$ into $E_M$. 

Suppose that $L>M$, i.e., that the map $\eta$ has a non-trivial
kernel, or that there exists a non-elliptic element in $L$ that is mapped to an elliptic element in $M$ by
$\eta$.

Let $M \to M_1 \to \ldots \to M_s$ be a (well-structured) resolution of $M$, i.e., 
a resolution of $M$ that satisfies the 
properties of the resolution that is associated with a limit group over free products in theorem 18. 
Then there exists a f.g.\ limit quotient of $L$, $CM$, with a set of elliptics, $E_{CM}$, and
a well-structured resolution of $CM$, $CM \to CM_1 \to \ldots \to CM_s$, 
that satisfies the properties of the resolutions
in theorem 18, 
and  a quotient map:
$\tau:L \to CM$, that maps $E_L$ into $E_{CM}$, so that:
\roster
\item"{(1)}" there exists a quotient map: $\nu: CM \to M$, that maps $E_{CM}$ onto $E_M$,
so that $\eta=\nu \circ \tau$.

\item"{(2)}"  
if $\eta: L \to M$ has a non-trivial
kernel, then $\tau:L \to CM$ has a non-trivial kernel. If there exists a non-elliptic element in $L$ that is mapped 
to an elliptic element in $M$ by
$\eta$, then there exists a non-elliptic element in $L$ that is mapped to an elliptic element in $CM$ by $\tau$.
If $M_{i+1}$ is a proper quotient of $M_i$, then $CM_{i+1}$ is a proper quotient of 
$CM_i$.

\item"{(3)}" if $\eta$ maps an elliptic element in $L$ to the identity, then $\tau$ maps an elliptic element in $L$
to the identity.

\item"{(4)}" if $M$ is a free product of finitely many elliptic subgroups and a free group, so is $CM$. More
generally, $CM_j$ is mapped onto $M_j$, $1 \leq j \leq s$, where elliptics in $CM_j$ are mapped onto
elliptics in $M_j$.

\item"{(5)}" all the homomorphisms of the given limit group $L$ that factor through the given well-structured
 resolution
of $M$, factor through the  resolution $CM \to CM_1 \to \ldots \to CM_s$.

\item"{(6)}" with the given well-structured resolution, $M \to M_1 \to \ldots \to M_s$, we can
naturally associate a completion, $Comp_M$ (see definition 1.12 in [Se2]), and with the resolution
$CM \to CM_1 \to \ldots \to CM_s$ we can naturally associate a completion, $Comp_{CM}$. $CM$ is embedded into
$Comp_{CM}$, and the elliptic elements in $CM$ are mapped into the terminal limit group $CM_s$. 

\item"{(7)}" by theorem 18, 
the elliptic elements, $E_L$,  in the limit group $L$ are conjugates of finitely many (possibly none) f.g.\ subgroups,
$E_1,\ldots,E_r$ in $L$. Then the completion, $Comp_{CM}$, is obtained from (copies of the subgroups) $E_1,\ldots,E_r$
by adding finitely many generators and relations, i.e., $Comp_{CM}$ is f.p.\ relative to the subgroups 
$E_1,\ldots,E_r$.  

\item"{(8)}" if $M$ admits a free decomposition, $M=V_1*\ldots*V_u*F_t$, where $F_t$ is a free group, and
this free decomposition is respected by the given resolution of $M$,
then $CM$ has a similar free decomposition, 
$CM=CV_1*\ldots*CV_u*F_t$, which is respected by the constructed  resolution of $CM$,
where the map $\nu$ respects this free decomposition, i.e., $\nu(CV_i)=V_i$,
$i=1,\ldots,u$, and $\nu(F_t)=F_t$. In particular, the completion, $Comp_{CM}$, admits a similar free
decomposition, $Comp_{CM}=Comp_1*\ldots*Comp_u*F_t$, where $CV_i$ embeds into $Comp_i$.
\endroster
\endproclaim

\nfp  Let $L$ be a limit group over free products, with set of elliptics $E_L$. By theorem 18 there are finitely
many subgroups, $E_1,\ldots,E_r$, in $L$, so that the set of elliptic elements in $L$, $E_L$, is the union of the
conjugacy classes of $E_1,\ldots,E_r$.  Let $M$ be a limit quotient of $L$, and let $M \to M_1 \to \ldots \to M_s$ be
a well-structured resolution of $M$, where $M_s$ is a free product of finitely many elliptic factors and a possibly
trivial free group. 

With the given well-structured resolution of $M$ we associate a completion, $Comp_M$. $M$ is a limit quotient of $L$,
and $M$ is a subgroup of the completion, $Comp_M$, so $L$ is mapped into $Comp_M$. Hence, the elliptic subgroups
in $L$, $E_1,\ldots,E_r$, are mapped into conjugates of the elliptic subgroups, that are factors in the terminal
limit group $M_s$, in $Comp_M$.  
If the terminal limit group $M_s$ is f.p.\ relative to the subgroups, $E_1,\ldots,E_r$, the theorem follows
(by taking the cover $CM$ to be $M$ and $Comp_{CM}$ to be $Comp_M$). Otherwise we modify the argument that was 
used 
to prove proposition 19.

Since $M$ is embedded into the completion $Comp_M$, $L$ is naturally mapped into $Comp_M$. 
Each of the factors of the terminal limit group of $Comp_M$, $M_s$,  is f.g.\ so it is a quotient of some conjugates
of (copies of) the elliptic subgroups of $L$, $E_1,\ldots,E_r$ and a f.g.\ free group. 
We start the construction of a the completion $Comp$ that covers the completion $Comp_M$,  
with a tower $T_1$ that has in 
its base level the free product of a free group (isomorphic to the free factor in the free decomposition 
of $M_s$), the free products of corresponding conjugates of $E_1,\ldots,E_r$ with
free groups (so that each of the factors of $M_s$ is a quotient of each of these free products). The next 
(upper) levels are connected to the lower levels
of the constructed tower $T_1$,
precisely as they are connected in the completion, $Comp_M$, i.e., using the same graphs of groups, just that the 
 group that is associated with the lowest level in $Comp_M$, which is $M_s$, is replaced by the prescribed free
products. 

$T_1$ is a tower, but it is not necessarily a completion (see definition 1.12 in [Se2]), 
as in general there are no retractions
from a group that is associated with a certain level onto the group that is associated with the level below it.
Each of the levels above the base level in $T_1$ is  constructed using a (finite) graph of groups, in which
some vertex groups are the groups that are associated with the lower level in $T_1$. Hence, the group that is
associated with a level above the base level,  
is obtained from a free product of the group that is associated with the lower level with a f.p.\ group
by imposing finitely many relations. Furthermore, the graphs of groups that are associated with the different
levels in $T_1$ are similar to the  graphs of groups that are associated with the corresponding levels in the
completion $Comp_M$, and differ from $Comp_M$ only in the groups that are associated with the base level. 

Each of the groups that are associated with the upper levels in $T_1$ is obtained from the groups
that appear in the lower level of $T_1$ by a free product with a f.g.\ free group and further 
imposing finitely many relations. The graphs of groups that are
associated with the upper levels in $T_1$ have similar structure as the corresponding graphs of groups that
are associated with the levels of the completion $Comp_M$,  i.e., the graphs of groups differ only
in the vertex groups that are associated with lower levels. Furthermore, these vertex groups differ 
only in the groups that are associated with the base levels of $T_1$ and $Comp_M$. Hence, 
it is enough to impose only finitely many (additional) 
relations 
from the defining relations of the various factors of the limit groups that is associated with
the base level of $Comp_M$, $M_s$, on the subgroup that is associated with the base level of $T_1$. 
This means imposing  finitely many (additional) relations on the associated free 
products of  free groups and conjugates of 
(copies of) the subgroups, $E_1,\ldots,E_r$, that form the group which is associated with the base level of $T_1$,
so that if we replace the group that
is associated with the base level of $T_1$,
 with the obtained  quotient, 
and construct from the obtained base subgroup a tower,
$T_2$, by imitating the construction of $Comp_M$ and $T_1$ (i.e., with similar graphs of groups in all the upper
levels), $T_2$ will be a completion (i.e., it is a tower with retractions between consecutive levels).
     
$T_2$ is a completion, but it may be that the limit  group $L$ is not mapped into it. 
$L$ is mapped into the completion $Comp_M$, and as a limit group it is finitely presented relative to its 
elliptic subgroups. Hence, once again, it is enough to impose 
only finitely many relations from the defining relations of the various factors of $M_s$, 
so that if we replace the group that is associated with the base level in $T_2$
with the obtained quotient,  and construct from it a tower $T_3$ by imitating the 
construction of the towers $Comp_M$, $T_1$, and $T_2$, $T_3$ is a completion, it is f.p.\ relative
to the elliptic subgroups, $E_1,\ldots,E_r$, and $L$ maps into it. 

We denote the images of the limit
group $L$ into the various levels of the completion $T_3$, by $CM, CM_1,\ldots,CM_s$.
By imposing finitely many additional relations on the base subgroup of $T_3$ from the relations of the
base subgroup of $Comp_M$, $M_s$, one can further guarantee that if $M$ is a proper quotient of $L$, then
$CM$ is a proper quotient of $L$, if $L>M$ then $L>CM$, and similarly, if $M_{j+1}$ is a proper
quotient of $M_j$ then $CM_{j+1}$ is a proper quotient of $CM_j$, and if $M_{j+1}>M_j$ then
$CM_{j+1}>CM_j$. We denote the obtained completion by $Comp_{CM}$ and its associated resolution as the resolution
that is associated with $CM$
(the obtained resolution is a well-structured resolution by construction).  
All the other properties of the  limit groups, and the associated  resolution and  completion, $Comp_{CM}$, that
are listed in the statement of the theorem follow easily from the construction.

\line{\hss$\qed$}

Given a f.g.\ limit group over free products $L$, and its limit quotient $M$ with an associated well-structured
resolution, $M \to M_1 \to \ldots \to M_s$, that satisfy 
the assumptions of theorem 24, and for which
$L>M$, we call a limit quotient $CM$ of $L$, 
that satisfies the conclusion of the theorem, a $cover$
of the limit quotient $M$, its associated well-structured resolution, $CM \to CM_1 \to \ldots \to CM_s$, a
$cover$ of the given resolution of $M$, and the associated completion, $Comp$, into which $CM$ is embedded, that was
constructed from the given well-structured resolution of $M$, a $cover$ $completion$. 

In constructing the Makanin-Razborov diagrams of a f.p.\ or a f.g.\ group over a free or a hyperbolic group,
we were able to show that the set of shortening quotients of a limit group over these groups
contain finitely many equivalence classes of maximal shortening quotients. In studying limit groups
over free products we are not able to prove a similar theorem. Over free products we
prove that given a limit group $L$, and fixing a cover for each pair of a   shortening quotient and its 
associated well-structured resolution, there exists
a finite subcollection of covers which is good for all the shortening quotients of $L$. 
As we will see in the sequel,
a similar statement on the existence of a finite subcollection of  cover completions (with a similar proof)
is sufficient for the construction of the Makanin-Razborov diagram over free products. 

\vglue 1pc
\proclaim{Theorem 25} Let $L$ be a f.g.\ limit group over free products, suppose that $L$ is not (entirely)
elliptic and that $L$ admits no free product decomposition in which the elliptic elements in $L$, $E_L$, can
be conjugated into the factors. 

With each pair of a   shortening quotient $M$ of $L$, and a well-structured resolution of $M$,
there is an associated quotient map, $\eta_M: L \to M$, 
that satisfies the assumptions of theorem 24. Hence, by the conclusion of theorem 24, for each 
pair of a shortening quotient $M$ of $L$, and its associated well-structured resolution we can choose a cover $CM(M)$
together with a completion, $Comp_{CM}$, into which $CM$ embeds.

From the entire collection of covers of  shortening quotients of $L$ and their associated well-structured resolutions,
 it is possible to choose
a finite subcollection of covers, $CM_1,\ldots,CM_e$,  so that for every maximal shortening quotient $M$,
there exists an index $i$, $1 \leq i \leq e$, for which the quotient map, $\eta: L \to M$, is a composition of
the two quotient maps: $L \to CM_i \to M$ (where elliptics are mapped to elliptics in these two maps).
\endproclaim

\nfp  The argument that we use is similar to the proof of the finiteness of the number
of equivalence classes of maximal shortening quotients (over free products) of a
f.p.\ group. Let $L$ be a f.g.\ limit group. By theorem 24, given a shortening
quotient of it, $M$, and a well structured resolution of that shortening quotient,
there exists a cover $CM$ of $M$, and $CM$ can be embedded into a completion, $Comp$,
that is obtained from the finitely many (conjugates of) elliptic subgroups of $L$,
by adding finitely many generators and relations. Therefore, there are at most countably
many such completions, $Comp$, and hence, at most countably many such covers, $CM$.

\noindent
Note that by lemma 22 each shortening quotient $M$ of $L$ is either entirely elliptic, 
or it is a proper quotient of $L$. In case the shortening quotient $M$ is not entirely elliptic, it follows
by theorem 24, that the associated cover, $CM$ is a proper quotient of $L$ (like the shortening quotient $M$). 

Suppose that the countable collection of covers does not contain a finite subcover, i.e.,
there is no finite subcollection of the constructed covers, $CM_1,\ldots,CM_e$, so that
for every shortening quotient $M$, the quotient map $L \to M$ factors as a composition of
quotient maps of limit groups over free products: $L \to CM_i \to M$,
for some index $i$, $1 \leq i \leq e$.

To contradict the lack of a finite subcover, we start by ordering the collection of covers,
$CM_1,CM_2,\ldots$. Since there is no finite subcover, there must be a sequence of indices,
$i_1,i_2,\ldots$, so that a shortening quotient, $M_{i_j}$, which is covered by $CM_{i_j}$,
is not covered by any of the
previous covers, $CM_1,\ldots,CM_{i_j-1}$.

\noindent
For each index $j$, the shortening quotient $M_{i_j}$ is a limit of shortest homomorphisms, and it is
not covered by any of the covers, $CM_1,\ldots,CM_{i_j-1}$. Hence, for each index $j$, there exists a
shortest homomorphism $h_j: L \to A_j*B_j$, 
that does not factor through any of the covers, $CM_1,\ldots,CM_{i_j-1}$.

We look at the sequence of homomorphisms $\{h_j\}$. A subsequence of this sequence converges into a limit
group (over free products) $R$, which is a quotient of the limit group $L$. Unless the limit group $R$ is
the (possibly trivial) free product of finitely many elliptic factors and a (possibly trivial) free group, 
a subsequence of the shortenings of these homomorphisms converges into
a shortening quotient $R_1$ of $R$, where the elliptics in $R$ are mapped monomorphically into the elliptics
in $R_1$, and $R_1$ is a proper quotient of $R$. By continuing iteratively and applying the d.c.c.\
for decreasing sequences of limit groups over free products (theorem 13),  we obtain a finite (strict) resolution
$R \to R_1 \to \ldots \to R_s$, where $R_s$ is a free product of finitely many f.g.\ elliptic subgroups and
a (possibly trivial) free group. For brevity,
we still denote the obtained subsequence of shortest homomorphisms, $\{h_j\}$.

The pair of the shortening quotient $R$, and its (strict) resolution, $R \to R_1 \to \ldots \to R_s$, is
one of the pairs of a shortening quotient of the limit group $L$, and its associated strict resolutions, with
which we have associated the covers, $CM_1,CM_2,\ldots$. Hence, one of these covers, $CM_r$, is a cover that is 
associated with this pair. Since a cover completion is finitely presented relative to the
elliptic subgroups of $L$, for large enough indices $j$, the homomorphisms $\{h_j\}$ factor 
through cover completion and hence factor through the cover $CM_r$. 
That contradicts the choice of the homomorphisms $\{h_j\}$, as for large $j$,
$h_j$ is  supposed not to factor through the covers, $CM_1,\ldots,CM_{j-1}$.

\line{\hss$\qed$}

Theorem 25 proves that given a f.g.\ limit group $L$ that admits no free decomposition in which the
elliptic elements, $E_L$, can be conjugated into the factors, it is possible to find finitely many  limit
quotients of $L$, one which is isomorphic to $L$ and is entirely elliptic, and the rest which are proper quotients
of $L$, that cover all its shortening quotients. This finite collection of covers is not canonical,
but in principle it can be taken as the next step in the Makanin-Razborov diagram. Except for the
entirely elliptic cover that is isomorphic to $L$, the other covers that are
associated with $L$  are all proper quotients of it, hence, in principle we
can continue with the construction iteratively. However, the d.c.c.\ that we proved is valid only for
sequences of strictly decreasing limit quotients, for which the quotients are proper and 
are monomorphic when restricted to 
elliptic elements (theorem 13). 

Therefore, to complete the construction of the Makanin-Razborov diagram of a f.p.\ group over free products
we use a different approach. Instead of constructing a finite cover of all the shortening quotients of a
given limit group (over free products), we construct a finite cover for all the (strict) resolutions of
the given limit group. With each strict resolution of the given limit group we associate a cover of that
resolution (which is a resolution by itself), and there are only countably many such covers, 
as the completion that is associated with the cover resolution is f.p.\ relative to the elliptic subgroups
of the original limit group. Then we use a similar argument to the one that was used in proving theorem
25 to prove that there exists a finite subcollection of the collection of cover resolutions, i.e., that there
exists a finite subcollection so that every homomorphisms of the given limit group into free products factors through
at least one of the resolutions from the finite subcollection of cover resolutions.

\vglue 1pc
\proclaim{Theorem 26} Let $L$ be a f.g.\ limit group over free products. Then there exists finitely many well-structured
 resolutions 
of quotients of $L$, so that every homomorphism from $L$ into a free product factors through at least one of 
these well-structured
resolutions. Furthermore, with each of these (finitely many) 
well-structured resolutions we can naturally associate a completion, and these completions are f.p.\ relative to the 
(finitely many) elliptic subgroups in the given limit group $L$.  
\endproclaim

\nfp  The proof is similar to the proof of theorem 25.
First, we factor the limit group over free products $L$ into a maximal free decomposition in which
the elliptic elements of $L$, $E_L$, can be conjugated into the factors. We continue with each of the factors
separately. Hence, we may assume that the limit group $L$ is freely indecomposable with respect to its set of elliptics,
$E_L$.
By theorem 24, given a limit quotient of $L$, that we denote $T$, and a well-structured resolution of $T$,
$T \to T_1 \to \ldots \to T_s$, that is obtained by taking successive shortening quotients (see theorem 18),
there exists a cover of $T$, that we denote $CT$, which is a limit quotient of $L$, 
and a cover of the resolution that is associated with $T$, which is a well-structured resolution, with
an associated cover completion, $Comp_{CT}$,
that is obtained from the finitely many (conjugates of) elliptic subgroups of $L$,
by adding finitely many generators and relations. Therefore, there are at most countably
many such triples of a cover of a limit quotient, an associated (well-structured) cover resolution, and the
corresponding cover completion.

Suppose that the countable collection of cover resolutions  does not contain a finite subcover, i.e.,
there is no finite subcollection of the constructed covers, $CT_1,\ldots,CT_e$, with associated
cover completions, $Comp_1,\ldots,Comp_e$, so that for each homomorphism $h$ of $L$ into a free product
(that maps the elliptics in $L$, $E_L$, into elliptic elements), the homomorphism $h$ factors through at least
one of the cover resolutions that is associated with the cover completions, $Comp_1,\ldots,Comp_e$. 

To obtain a contradiction to the lack of finiteness of covering resolutions,  we start by ordering the 
collection of covering completions and their associated resolutions, $Comp_1,Comp_2,\ldots$.
Since there is no finite subcover for the entire collection of homomorphisms of the given limit group $L$
into free products, for each index $i$, there exists a homomorphism, $h_i:L \to A_i*B_i$, that does not factor
through the resolutions that are associated with the completions, $Comp_1,\ldots,Comp_{i-1}$.

Like in the proof of theorem 25, a subsequence of the sequence of homomorphisms, $\{h_i\}$, 
converges into a limit
group (over free products) $R$, which is a quotient of the limit group $L$. Unless $R$ is a (possibly trivial)
free product
of elliptic subgroups and a (possibly trivial) free group, 
a subsequence of the shortenings of these homomorphisms converges into
a shortening quotient $R_1$ of $R$, where the elliptics in $R$ are mapped monomorphically into the elliptics
in $R_1$, and $R_1$ is a proper quotient of $R$. By continuing iteratively and applying the d.c.c.\
for decreasing sequences of limit groups over free products (theorem 13),  we obtain a finite well-structured
 resolution
$R \to R_1 \to \ldots \to R_s$, where $R_s$ is a free product of finitely many f.g.\ elliptic subgroups and
a (possibly trivial) free group. For brevity,
we still denote the obtained subsequence of shortest homomorphisms, $\{h_i\}$.

The pair of the limit quotient $R$ of the given limit group (over free products) $L$, 
and its (well-structured) resolution, $R \to R_1 \to \ldots \to R_s$, is
one of the pairs of a limit quotient of  $L$, and its associated well-structured resolutions, with
which we have associated the covers, $Comp_1,Comp_2,\ldots$. Hence, one of these completions, $Comp_r$, 
is a cover that is 
associated with this pair. Therefore, for large enough index $i$, the homomorphism $\{h_i\}$ factors 
through the cover resolution that is associated with the completion, $Comp_r$. 
That contradicts the choice of the homomorphisms $\{h_i\}$, as for each $i$,
$h_i$ is  supposed not to factor through the cover resolutions that are associated with the completions,
$Comp_1,\ldots,Comp_{i-1}$.

\line{\hss$\qed$}

\vglue 1.5pc
\centerline{\bf{\S5. Makanin-Razborov Diagrams of Finitely Presented Groups}}
\medskip

Theorem 21 on the finiteness of the number of equivalence classes of maximal limit quotients 
(over free products) of a f.p.\
group, together with theorem 26 on the existence of finitely many (cover) resolutions of some quotients of
a given f.g.\ limit group over free products, so that every homomorphism of the given f.g.\ limit group
into free products factors through at least one of the resolutions, allow us to construct a Makanin-Razborov
diagram of a f.p.\ group over free products.

Given a f.p.\ group $G$, we start with its (canonical) finite collection of maximal limit quotients
over free products (theorem 21). With each maximal limit quotient we associate a finite collection of 
well-structured resolutions of it (according to theorem 26), so that each homomorphism of the original 
maximal limit quotient into free products, factors through at least one of its associated resolutions. We
construct the diagram by mapping the given f.p.\ group $G$ into the f.g.\ limit group that appears in the
top level of each of the (finitely many) well-structured resolutions that are associated with its collection
of maximal limit quotients (in parallel). Since every homomorphism of $G$ into free products, factors through
at least one of its maximal limit quotients, every homomorphism of $G$ into free products 
factors through at least one of the 
resolutions in its Makanin-Razborov diagram over free products. That is for every homomorphism 
of the f.p.\ group $G$, there exists at least one resolution in the Makanin-Razborov diagram, so that
the homomorphism can be written 
as a successive composition of the epimorphisms between the groups that appear in the various
levels of the resolutions, modular automorphisms of the limit groups that appear in the various levels
(that are encoded by the virtually abelian decompositions that are associated with these groups), and finally
a homomorphism from the terminal group of the resolution (which is a free product of elliptic factors and
a free group), that sends every elliptic factor into a conjugate of a factor in the image free product.  

At this stage we slightly improve the diagram. The virtually abelian decompositions that are associated with each of the limit groups that appear in the various levels of the well-structured resolutions in the
Makanin-Razborov diagrams, are decompositions that are inherited from the free and virtually abelian 
JSJ decompositions of the limit groups that appear along the well-structured resolutions that 
the resolutions in the Makanin-Razborov diagram  cover, according to the construction that appears in
theorem 24. However, these may not be the Grushko and virtually abelian decompositions of the limit
groups in the Makanin-Razborov diagram themselves. To fix that, and make sure that all the decompositions
in the limit groups that appear in the Makanin-Razborov diagram are indeed Grushko and virtually
abelian JSJ decompositions (over free products), we slightly modify the construction of a cover.

\vglue 1pc
\proclaim{Theorem 27} Let $L$ be a f.g.\ limit group over free products, let $M$ be a limit quotient
of $L$, and let $M \to M_1 \to \ldots \to M_s$, be a well-structured resolution of $M$, so that $M_s$
is a free product of finitely many elliptic factors and a possibly trivial free group. Suppose that the free
products that are associated with the various limit groups along the resolution, $M,M_1,\ldots,M_s$, are
their Grushko free decompositions with respect to their elliptic subgroups (i.e.,  the
resolution respects the Grushko free decompositions of the groups along it), 
and that the virtually abelian decompositions that
are associated with the limit groups $M,M_1,\ldots,M_s$ are their virtually abelian JSJ decompositions over
free products.

Then there exists a cover $CM$ of $M$, with a cover resolution, $CM \to CM_1 \to \ldots \to CM_s$,
that satisfies the properties of a cover that are listed in theorem 24, and for which the free decompositions
along the cover resolution are the Grushko free decompositions of the limit groups, $CM,CM_1,\ldots,CM_s$,
and the virtually abelian JSJ decompositions of these groups over free products have the same
structure as the virtually abelian decompositions that are associated with them along the resolution, i.e., the same
structure as the virtually abelian JSJ decompositions of the limit groups, $M,M_1,\ldots,M_s$.  
\endproclaim

\nfp The proof that we use is a modification of the argument that we used to prove theorem 24.
 Let $L$ be a limit group over free products, with a set of elliptics $E_L$. Recall that by theorem 18, the set
of elliptics $E_L$ is the union of conjugates of some (elliptic) subgroups, $E_1,\ldots,E_r$, in $L$.
Let $M$ be a limit quotient of $L$, and let $M \to M_1 \to \ldots \to M_s$ be
a well-structured resolution of $M$, where $M_s$ is a free product of finitely many elliptic factors and a possibly
trivial free group. 

With the given well-structured resolution of $M$ we associate a completion, $Comp_M$. Given the well-structured resolution
of $M$, and its associated completion, $Comp_M$, we use the construction that was used in proving theorem 24, and
construct a completion, $Comp$, which is f.p.\ relative to the elliptic subgroups, $E_1,\ldots,E_r$, and for which
the images of the limit group $L$ into the various levels of $Comp$, that were denoted, $CM, CM_1,\ldots,CM_s$,
 satisfy the list of properties that is
presented in theorem 24.

By adding finitely many relations to the base subgroup of $Comp$ from the set of relations that are
defined on the base subgroup, $M_s$, of the completion $Comp_M$, we may assume that the abelian decompositions
that are inherited by the subgroups, $CM,CM_1,\ldots,CM_s$, from the abelian
decompositions that are associated with the various levels of the completion $Comp$,
are similar to the abelian decompositions that are inherited by the various abelian decompositions of the
subgroups, $M_1,\ldots,M_s$ from the abelian decompositions that are associated with the various levels of 
$Comp_M$.

Suppose that the Grushko free decomposition of the limit group $M$ with respect to its elliptic subgroups is
$M=M^1*\ldots*M^b*F_v$, and this free decomposition together with the virtually abelian JSJ decompositions of
the factors, $M^j$, over free products with respect to the elliptic subgroups of $M$, 
give rise to an abelian decomposition , $\Delta_M$. Note that by our assumptions, the completion, $Comp_M$,
 respects the Grushko free decomposition of $M$, and the abelian decompositions that are associated
with the various levels of $Comp_M$ are the virtually abelian JSJ decompositions over free products of the
subgroups, $M,M_1,\ldots,M_{s-1}$.  

We order the relations that the terminal limit group $M_s$ of $Comp_M$ satisfy, and sequentially impose them on
the terminal limit group of the completion, $Comp$. We claim that after adding finitely many of these relations,
the free product 
decomposition, and the virtually abelian JSJ decomposition of the corresponding subgroup $CM$ (after adding
the relations) will be similar to those of the subgroup $M$.  

The cover $CM$, which is the image of the limit group $L$ in the completion, $Comp$, admits a free 
decomposition $CM=CM^1*\ldots*CM^b*F_v$, in which the elliptic subgroups in $CM$ can be conjugated into
the factors.  This free decomposition  is inherited from the structure of the completion, $Comp$,
as the completions $Comp$ and $Comp_M$ have the same structure, and $Comp_M$ respects the Grushko decomposition
(relative to elliptic subgroups) of the limit quotient $M$, 
$M=M^1*\ldots*M^b*F_v$. 

Let $CM(n)$ be the image of $L$ in the completion, $Comp(n)$, that is obtained from $Comp$ 
by imposing on the terminal
level in $Comp$ the first $n$ relations in $M_s$, the terminal limit group in $Comp_M$. $CM(n)$ inherits 
a free decomposition from $Comp(n)$, 
$CM(n)=CM(n)^1*\ldots*CM(n)^b*F_v$, a free decomposition in which the elliptic subgroups in $CM(n)$ can be 
conjugated into
the factors  (note that the elliptic subgroups in $CM(n)$ can be conjugated into the factors of the terminal 
limit group of $Comp(n)$). If this free decomposition is not the Grushko free decomposition of $CM(n)$ with
respect to its elliptic subgroups, 
then at least one of the factors admits a further  non-trivial  free decomposition with respect to the
elliptic subgroups.  

Suppose that there exists a sequence of indices (still denoted $n$) for which the free decomposition of $CM(n)$
that is inherited from $Comp(n)$ is not the Grushko free decomposition of $CM(n)$ with respect to the elliptic
subgroups in $CM(n)$. By passing to a subsequence (still denoted $n$) 
we may assume that one of the factors, w.l.o.g. $CM(n)^1$
admits a non-trivial free decomposition $CM(n)^1=A_n*B_n$, where each of the elliptic subgroups in $CM(n)$ can
be conjugated into one of the other factors in the given free decomposition of $CM(n)$, to $A_n$ or to $B_n$.

In that case we look at the actions of the groups $CM(n)^1$ on the (pointed) Bass-Serre trees,
$(T_n,t_n)$, that correspond to the 
(non-trivial) free products, $A_n*B_n$. Note that these actions are faithful actions of the groups, $CM(n)^1$, 
that the elliptic subgroups in $CM(n)$ that can be conjugated into $CM(n)^1$ can be conjugated into $A_n$ 
or $B_n$,
and that by construction, the
direct limit of the groups, $CM(n)^1$, is the factor $M^1$ of the limit group $M$ which is assumed to be
freely indecomposable
relative to its elliptic subgroups.

$CM^1$ is f.g.\ so we fix a generating set for it, $<g_1,\ldots,g_d>$,
and since the groups $CM(n)^1$ are (limit) quotients of
$CM^1$, it gives us a generating set for each of the groups, $CM(n)^1$.  Given the  action of $CM(n)^1$ on the
Bass-Serre tree, $(T_n,t_n)$, we precompose this action with a (modular) automorphism $\phi_n$ of $CM(n)^1$, i.e., an 
automorphism that can be expressed as a composition of an automorphism that comes from the virtually abelian
decomposition that $CM(n)^1$ inherits from the virtually abelian decomposition that is  associated
with the top level in $Comp(n)$ and an inner automorphism,  so that the maximal displacement of the base point
$t_n$ by the action of the tuple of elements, $\phi_n(g_1),\ldots,\phi_n(g_d)$, is minimal among all such 
(modular) automorphisms $\phi$.

Since we modify the actions of the groups, $CM(n)^1$, by precomposing them with (modular) automorphisms, 
and since the actions are all faithful,
there is a subsequence of twisted actions that converge into an action of the direct limit of the
groups, $CM(n)^1$, i.e., the factor $M^1$ of $M$, on a real tree. Since the automorphisms $\phi_n$ were chosen
to minimize the displacement of the base points under the corresponding twisted actions, and since the virtual
abelian JSJ decomposition of the limit group $M^1$ has the same structure as the virtually abelian 
decomposition that is inherited by $CM(n)^1$ from the virtually abelian decomposition that is
associated with the top level of the completions, $Comp(n)$, the set of displacements
of the base points under the twisted actions has to be bounded. Hence, the factor $M^1$ of $M$ inherits a non-trivial free
decomposition from the limit action, a free decomposition in which all the elliptic subgroups in $M^1$
can be conjugated into the factors. This contradicts the assumption that $M^1$ admits no such non-trivial
free decomposition. Therefore, there must exist some index $n_0$, so that for all $n>n_0$, the limit groups
$CM(n)^1$ admit no free decomposition in which the elliptic subgroups of $CM(n)^1$ can be conjugated into the 
factors. 

By passing to a subsequence, we may assume that all the factors in the free decomposition of the limit groups,
$CM(n)$, are freely indecomposable relative to their elliptic subgroups.
Suppose that there exists a sequence of indices (still denoted $n$) for which the virtually abelian
decomposition that at least one of the factors of the the groups, $CM(n)$, $CM(n)^1,\ldots,CM(n)^b$,
inherits from the virtually abelian decomposition that is associated with the top level
of the completion, $Comp(n)$, is not the virtually abelian JSJ decomposition over free products of that factor.
Wlog we may assume that this factor is $CM(n)^1$. 

Let $\Delta(n)$ be the
virtually abelian decomposition that $CM(n)^1$ inherits from the virtually abelian 
decomposition that is associated with
the top level of the completion $Comp(n)$. Let $JSJ(n)$ be the virtually abelian JSJ decomposition 
over free products of
$CM(n)^1$, and let $\Delta_M$ be the virtually abelian decomposition that $M$ inherits from the virtually abelian
decomposition that is associated with the top level of the completion, $Comp_M$, which by our assumptions is
the virtually abelian JSJ decomposition of $M$ over free products. Since we assumed that the virtually abelian
decompositions, $\Delta(n)$, are not identical to the virtually abelian decompositions, $JSJ(n)$, the virtually
abelian JSJ decompositions, $JSJ(n)$, must be proper refinements of the virtually abelian decompositions,
$\Delta(n)$. Note that the structure of the virtually abelian decompositions, $\Delta(n)$, is similar to that
of the abelian decomposition, $\Delta_M$.  

For every index $n$, the virtually abelian JSJ decomposition $JSJ(n)$ is a proper refinement of the 
virtually abelian  decomposition $\Delta(n)$. Hence, if needed we can cut some of the $QH$ subgroups in 
$JSJ(n)$ along s.c.c.\ and obtain a new decomposition, $\Theta(n)$, of $CM(n)^1$ that refines $\Delta(n)$,
in which all the edge groups and all the $QH$ vertex groups in $\Delta(n)$ are elliptic, and at least
one of the non-$QH$ non-virtually-abelian vertex groups in $\Delta(n)$ is not elliptic. 
Hence, at least one of these 
vertex groups inherits a non-trivial virtually abelian decomposition from $\Theta(n)$, a decomposition in
which all the edge groups that are connected to that vertex group are elliptic.

\noindent
By passing to a further subsequence (still denoted $n$), we may assume that the vertex group that inherits a 
non-trivial virtually abelian decomposition from $\Theta(n)$ is a vertex group $V(n)$ in $\Delta(n)$
that is mapped to the same vertex 
group $V$ in $\Delta_M$, the virtually abelian JSJ decomposition of the limit group $M$. 

We  fix a free group $F_r$, where $r$ is
the rank of the limit group $CM$, and an epimorphism, $\tau: F_r \to CM$.
We fix a finite generating set for $F_r$. We may assume that this generating
set contains elements that are mapped to elements that generate the edge groups and the
vertex groups in the virtually 
abelian decomposition of $CM$ that is
inherited from the top level of the completion, $Comp$. 

For each index $n$, we
look at a homomorphism $h_n:F_r \to A_n*B_n$ that approximates the limit group $CM(n)$. This means
that  $h_n$ maps 
each element in the ball of radius $n$ in the Cayley graph of $F_r$ (with respect to the given set of generators), to
an elliptic element or to a trivial element if and only if the element is trivial or elliptic in $CM(n)$. 
It maps the elements from the generating sets that are mapped to the edge groups in $\Delta(n)$
to  non-elliptic elements. Furthermore, let $S<F_r$ be the subgroup that is generated by those elements
in the fixed generating set of $F_r$ whose image generate the vertex group in the virtually abelian
decomposition of $CM$ that is mapped to the vertex group $V$ in $\Delta_M$, and the edge groups that
are connected to that vertex group. The vertex group $V(n)$ is not elliptic in the virtually abelian 
decomposition of the factor $CM(n)^1$, $\Theta(n)$, and the edge groups that are connected to $V(n)$ in
$\Delta(n)$ are elliptic in $\Theta(n)$. Hence, we may further modify each of the
homomorphisms $h_n$, by precomposing each of them  with Dehn twists along edge groups that lie in the graph of groups
that is inherited by $V(n)$ from the graph of groups $\Theta(n)$. We apply this modification,
so that for the obtained homomorphism, $\hat h_n$, when restricted to the subgroup $S<F_r$ (which is mapped onto $V(n)$),
 $\hat h_n:S \to A_n*B_n$, the minimal displacement of a
point in the Bass-Serre tree, that is associated with the free product $A_n*B_n$, under the action of the
fixed set of generators
of $S$, will be at least $n$ times larger than the minimal displacement of a point in that Bass-Serre
tree, under the action of the fixed set
of generators that are mapped to any given edge group that is connected to $V(n)$ in $\Delta(n)$.

By construction, the homomorphisms, $\{\hat h_n:S \to A_n*B_n\}$, converge into a non-trivial
 action of the vertex group $V$
in $\Delta_M$ on some real tree (where the convergence is into $V$ as a limit group over free products). 
All the edge groups that are connected to $V$ in $\Delta_M$ fix
points in that real tree and they are all non-elliptic subgroups (i.e. each element in these groups is
mapped to non-elliptic element in $A_n*B_n$ for large $n$). With this action it is possible to associate
a non-trivial graph of groups decomposition of $V$, with abelian edge groups, in which all the edge groups
that are connected to $V$ are contained in vertex groups in that graph of groups decomposition. Hence,
using this graph of groups decomposition it is possible to further refine the graph of groups, $\Delta_M$,
and this clearly contradict the assumption that $\Delta_M$ is the virtually abelian JSJ decomposition of
the limit group $M$. 

Therefore, for large $n$, the abelian decompositions, $\Delta(n)$, are the virtually
abelian JSJ decompositions of the limit groups over free products, $CM(n)$. The same argument implies
the same results for the next limit groups in the constructed resolution, $CM_1(n),\ldots,CM_{s-1}(n)$,
and the theorem follows.

\line{\hss$\qed$}

The Makanin-Razborov diagram of a f.p.\ group $G$ over free products is uniform, i.e., it  encodes 
all the homomorphisms from
$G$ into arbitrary free products. Equivalently, it encodes all the quotients of a f.p.\ group that
are free products. As we will see in the sequel, the Makanin-Razborov diagram that we constructed
suffices 
in order to modify the results and the techniques that
were used to study the first order theory of a free or a hyperbolic group, in order to study the first order
theory of a free product. We also believe that modifications of it can be applied for studying homomorphisms
of a f.p.\ group into groups with more general splittings (notably $k$-acylindrical splittings), and probably
homomorphisms into (some classes of) relative hyperbolic groups.

Unfortunately, the diagram that we constructed is not canonical, as it uses finite covers (theorems
25 and 26), and these are not 
unique. To construct a canonical diagram, we believe that it's better to study only maximal homomorphisms into
free products.

\vglue 1pc
\proclaim{Definition  28} Let $G$ be a f.g.\ group.
On the set of  homomorphisms of $G$ into free products, we define a partial order.
Let $h_i:G \to A_i*B_i$, $i=1,2$,  be two homomorphisms. Note that the images of the homomorphisms $h_i$ inherit
(possibly trivial) free products from the free product decompositions $A_i*B_i$, $i=1,2$.
We write that  $h_1 > h_2$, 
if there exists an  epimorphism with non-trivial kernel: $\tau: h_1(G) \to h_2(G)$, 
that maps the elliptics in $h_1(G)$ into the elliptics
in $h_2(G)$,  
so that for every $g \in G$, $h_2 (g)= \tau (h_1(g))$. 

If $\tau$ is an isomorphism and it maps the elliptics in $h_1(G)$
onto the elliptics in $h_2(G)$, and for every $g \in G$, 
$h_2(g) = \tau ( h_1(g))$, 
we say that $h_1$ is in the same equivalence class as $h_2$.

Note that this relation on homomorphisms into free products, which is a partial order on
homomorphisms, is a special case
of the partial order that was defined in Definition 12 for limit groups over free products. 
\endproclaim

To construct a canonical Makanin-Razborov diagram, it seems that one needs to prove the existence of maximal
homomorphisms with respect to the above partial order. The existence of maximal homomorphisms allows one
to construct a canonical (finite) collection of maximal shortening quotients of a f.g.\ limit group over
free products, and then prove a d.c.c.\ that allows the termination of the construction of a diagram,
using somewhat similar construction to the one used over free and hyperbolic groups.
To prove the existence of maximal homomorphisms (with respect to the prescribed partial order), one
needs to prove the following natural conjecture: 

\vglue 1pc
\proclaim{Conjecture} Let $G$ be a f.g.\ group. 
Let $h_1,h_2,\ldots$ be a sequence of homomorphisms of $G$ into free products, for which:
 $$h_1 \, < \, h_2 \, < \, \ldots$$ 
 Then there exists a homomorphism $h$ from $G$ into a free product,
so that for every index $m$, $h > h_m$ (one may even assume that the homomorphisms, $\{h_m\}$, 
do not factor through
an epimorphism onto a group of the form $M*F$ for some nontrivial free group $F$). 
\endproclaim

Finally, we note that the Makanin-Razborov diagram over free products that we constructed is associated
with a f.p.\ group. Some of our arguments are not valid for f.g.\ groups. In particular, although there exist
maximal elements in the set of limit quotients over free products of a f.g.\ group, 
it is not clear if there are only finitely many
maximal limit quotients. 
Therefore, the study of the collection
of homomorphisms from a given f.g.\ group into free products remains open.


\vskip 2em

\end